\documentclass[a4paper,reqno]{amsart}
\usepackage{lmodern}
\usepackage{lmodern}
\usepackage[T1]{fontenc}
\usepackage[utf8]{inputenc}
\usepackage{xcolor}
\usepackage{mathrsfs}
\usepackage{mathtools}
\usepackage{amstext}
\usepackage{amsthm}
\usepackage{amssymb}
\usepackage{microtype}
\usepackage[bookmarks=false,
 breaklinks=false,pdfborder={0 0 1},backref=false,colorlinks=true]
 {hyperref}
\hypersetup{
 linkcolor=black,citecolor=black,urlcolor=black}

\makeatletter

 \pdfpagewidth\paperwidth

\numberwithin{equation}{section}
\numberwithin{figure}{section}

\usepackage{mathrsfs}
\usepackage{aliascnt}
\usepackage[capitalise,noabbrev]{cleveref}
\usepackage{pgfplots}
\usepgfplotslibrary{groupplots}
\pgfplotsset{compat=1.17}

\renewcommand{\phi}{\varphi}
\newcommand{\N}{\mathbb N}
\newcommand{\Z}{\mathbb Z}
\newcommand{\R}{\mathbb R}
\newcommand{\cA}{\mathcal A}

\newcommand{\cM}{\mathcal M}
\newcommand{\cP}{\mathcal P}
\newcommand{\Rep}{\mathfrak R}
\newcommand{\dd}{\,\mathrm d}
\newcommand{\e}{\,\mathrm e}

\theoremstyle{plain}
\newtheorem{thm}{Theorem}[section]
\newaliascnt{prop}{thm}
\newtheorem{prop}[prop]{Proposition}\aliascntresetthe{prop}
\newaliascnt{cor}{thm}
\newtheorem{cor}[cor]{Corollary}\aliascntresetthe{cor}
\newaliascnt{lem}{thm}
\newtheorem{lem}[lem]{Lemma}\aliascntresetthe{lem}
\theoremstyle{definition}
\newaliascnt{defn}{thm}
\aliascntresetthe{defn}
\newaliascnt{assumption}{thm}
\aliascntresetthe{assumption}
\newaliascnt{example}{thm}
\newtheorem{example}[example]{Example}\aliascntresetthe{example}
\theoremstyle{remark}
\newaliascnt{rem}{thm}
\newtheorem{rem}[rem]{Remark}\aliascntresetthe{rem}

\AddToHook{env/prop/begin}{\crefalias{thm}{prop}}
\AddToHook{env/lem/begin}{\crefalias{thm}{lem}}
\AddToHook{env/cor/begin}{\crefalias{thm}{cor}}
\AddToHook{env/rem/begin}{\crefalias{thm}{rem}}
\AddToHook{env/example/begin}{\crefalias{thm}{example}}
\AddToHook{env/assumption/begin}{\crefalias{thm}{assumption}}
\crefname{cor}{corollary}{corollaries}
\Crefname{cor}{Corollary}{Corollaries}
\crefname{example}{example}{examples}
\Crefname{example}{Example}{Examples}
\crefname{lem}{lemma}{lemmas}
\Crefname{lem}{Lemma}{Lemmas}
\crefname{prop}{proposition}{propositions}
\Crefname{prop}{Proposition}{Propositions}
\crefname{rem}{remark}{remarks}
\Crefname{rem}{Remark}{Remarks}
\crefname{thm}{theorem}{theorems}
\Crefname{thm}{Theorem}{Theorems}
\crefname{defn}{definition}{definitions}
\Crefname{defn}{Definition}{Definitions}
\crefname{assumption}{assumption}{assumptions}
\Crefname{assumption}{Assumption}{Assumptions}

\definecolor{lime}{HTML}{A6CE39}
\DeclareRobustCommand{\orcidicon}{%
	\begin{tikzpicture}
	\draw[lime, fill=lime] (0,0) 
	circle [radius=0.16] 
	node[white] {{\fontfamily{qag}\selectfont \tiny ID}};
	\draw[white, fill=white] (-0.0625,0.095) 
	circle [radius=0.007];
	\end{tikzpicture}
	\hspace{-2mm}
}

\foreach \x in {A, ..., Z}{%
	\expandafter\xdef\csname orcid\x\endcsname{\noexpand\href{https://orcid.org/\csname orcidauthor\x\endcsname}{\noexpand\orcidicon}}
}

\makeatother

\usepackage{cleveref}
\AddToHook{env/lem/begin}{\crefalias{thm}{lem}}
\AddToHook{env/prop/begin}{\crefalias{thm}{prop}}
\AddToHook{env/rem/begin}{\crefalias{thm}{rem}}
\AddToHook{env/cor/begin}{\crefalias{thm}{cor}}
\AddToHook{env/example/begin}{\crefalias{thm}{example}}
\crefname{cor}{corollary}{corollaries}
\Crefname{cor}{Corollary}{Corollaries}
\crefname{example}{example}{examples}
\Crefname{example}{Example}{Examples}
\crefname{lem}{lemma}{lemmas}
\Crefname{lem}{Lemma}{Lemmas}
\crefname{prop}{proposition}{propositions}
\Crefname{prop}{Proposition}{Propositions}
\crefname{rem}{remark}{remarks}
\Crefname{rem}{Remark}{Remarks}
\crefname{thm}{theorem}{theorems}
\Crefname{thm}{Theorem}{Theorems}

\begin{document}
\title[Dimension spectra for random walks with a reflecting boundary]{Dimension spectra of escaping sets for one-dimensional random walks
with a reflecting boundary}
\author{Johannes Jaerisch\orcidB{}}
\address{Graduate School of Mathematics, Nagoya University, Furocho, Chikusaku,
Nagoya, 464-8602 Japan}
\email{jaerisch@math.nagoya-u.ac.jp}
\author{Marc Kesseböhmer\orcidC{}}
\address{Institute for Dynamical Systems, Faculty 3 -- Mathematics and Computer
Science, University of Bremen, Bibliothekstr. 5, 28359 Bremen, Germany}
\email{mhk@uni-bremen.de}
\date{\today}
\thanks{JJ was supported by JSPS KAKENHI Grant Numbers 21K03269 and 24K06777.
MK acknowledges the support through project L3 of the Collaborative
Research Center TRR 181 \emph{``Energy Transfer in Atmosphere and
Ocean''}, funded by the German Research Foundation (DFG; grant agreement
no. 274762653).}
\begin{abstract}
We study the $\mathbb{N}$-extension of a finite-branch uniformly
expanding interval map, whose fibre coordinate evolves as a one-dimensional
random walk reflected at the origin by Lindley reflection. For a prescribed
escape rate $\alpha$, we compare the infinitely often (i.o.) $\alpha$-escaping
set, whose fibre coordinate returns infinitely often to a bounded
window around $\alpha n$, with the uniformly $\alpha$-escaping set,
whose fibre coordinate never leaves such a window. We show that the
Hausdorff dimension spectra of the i.o.\ and uniform escaping sets
do not necessarily coincide. The strict-gain set, where the i.o.\
spectrum is strictly larger than the uniform one, is a finite union
of intervals. The uniform spectrum is real analytic for positive escape
rates, while the i.o.\ spectrum is real analytic away from the boundary
points of the strict-gain set, where its phase transitions are characterized
explicitly.
\end{abstract}

\maketitle

\section{Introduction}

The fibre recursion is the classical Lindley recursion from queueing
theory~\cite{Lindley52}, for which standard background is provided
by~\cite{Asmussen03}. In the locally constant full-shift setting,
a Bernoulli measure on the base makes the steps independent and identically
distributed, so that the fibre coordinate is an ordinary reflected
random walk on $\N$; more general Gibbs measures replace the independent
steps by a stationary, exponentially mixing symbolic process. Classical
probabilistic results already show that return behaviour on the half-line
depends sharply on the drift. For random walks that are homogeneous
outside a finite reflecting zone, Lalley~\cite{Lalley95} obtained
the following trichotomy for fixed-state return probabilities: negative
drift gives positive recurrence, zero drift gives a decay of order
$n^{-1/2}$, and positive drift gives a decay of order $R^{-n}n^{-3/2}$
for some $R>1$. Related reflection mechanisms and fluctuation results
were subsequently studied in~\cite{EssifiPeigneRaschel13}; see~\cite{IksanovEtAl25}
for a recent account of locally perturbed and reflected random walks.
At the pathwise level, large deviations for Lindley processes are
likewise organized by long atypical excursions away from the boundary~\cite{DuffyMeyn10}.
This is closely related to the boundary-waiting--excursion--return
mechanism underlying the dimension gap studied below. Our purpose,
however, is different from that of probabilistic return and limit
theorems: we determine the Hausdorff dimensions of exceptional trajectories
with a prescribed linear escape rate, distinguishing between i.o,\
and uniform escape.

A reflecting boundary can change the dimension of an escaping set
even when the prescribed escape rate is positive. An orbit may wait
near the boundary, make a long positive excursion, and return before
the next excursion. Such paths contribute to infinitely often escape
but not to uniform escape within a bounded window. This paper determines
the two Hausdorff dimension spectra and identifies exactly when boundary
returns increase the dimension. While our previous work~\cite{GJK-reflected}
treated the zero escape rate, where the i.o.\ and uniform escaping
sets have the same Hausdorff dimension, the present paper reveals
a new effect in transient dynamics: reflection can create a genuine
gap between the two spectra.

Let $F_{0}:\Delta_{0}\to(0,1)$, with $\Delta_{0}:=\bigcup_{i\in I}J_{i}\subset(0,1)$,
be a uniformly expanding full-branch interval map with finitely many
pairwise disjoint, non-empty open branch intervals. Let $\Rep_{0}$
be the invariant set of points whose forward orbits remain in $\Delta_{0}$.
The regularity, coding and distortion assumptions are specified in
\Cref{sec:setting}.

Throughout, $\N=\{0,1,2,\ldots\}$, and for $x\in\Rep_{0}$, its branch
itinerary $(x_{n})_{n\in\N}$ is defined by $x_{n}=i$ whenever $F^{n}_{0}(x)\in J_{i}$.
Thus $I$ is also the symbolic alphabet, and $x_{0}$ denotes the
branch containing $x$.

Assign an integer step $\psi(i)$ to each branch $J_{i}$, and define
$\psi(x):=\psi(x_{0})$. In symbolic terminology, $\psi$ is constant
on one-cylinders: its value depends only on the first symbol of the
itinerary. Its finite step set is $A:=\{\psi(i):i\in I\}\subset\Z$.
We assume throughout that 
\[
a_{-}:=\min A<0<a_{+}:=\max A,\qquad\gcd\{|a|:a\in A\}=1.
\]
Thus both positive and negative steps occur, and the steps generate
$\Z$ as an additive group.

The main object of this paper is the $\N$-extension 
\[
F_{\N}:\Rep_{0}\times\N\longrightarrow\Rep_{0}\times\N,\qquad F_{\N}(x,k)=\bigl(F_{0}x,\max\{0,k+\psi(x)\}\bigr).
\]
The second coordinate records the height. An attempted step below
zero is replaced by a move to zero; this boundary rule is called Lindley
reflection. Equivalently, writing 
\[
\rho_{i}(k):=\max\{0,k+\psi(i)\},\qquad i\in I,\quad k\in\N,
\]
we have $F_{\N}(x,k)=(F_{0}x,\rho_{x_{0}}(k))$.

The full-branch structure and the assumptions on $A$ imply that the
$\N$-extension is topologically transitive; see \Cref{lem:transitive}.

Since the base is contained in $(0,1)$, we identify $(x,k)$ with
the real number $x+k$. We use the same notation $F_{\N}$ for the
resulting interval map, so that 
\[
F_{\N}(x+k)=F_{0}(x)+\rho_{x_{0}}(k),
\]
and write $h(x+k):=k$ for the fibre coordinate.

For each $\alpha\ge0$ we consider the\emph{ infinitely often (i.o.)
$\alpha$-escaping set
\[
E^{\N}(\alpha)\coloneqq\left\{ x:\,\,\liminf_{n\to\infty}|h(F^{n}_{\N}x)-\alpha n|<\infty\right\} 
\]
}and the \emph{uniform $\alpha$-escaping set
\[
E^{\N}_{u}(\alpha)\coloneqq\left\{ x:\,\,\limsup_{n\to\infty}|h(F^{n}_{\N}x)-\alpha n|<\infty\right\} .
\]
}Here and below $x$ ranges over $\Rep_{0}+\N$, the phase space of
the interval map $F_{\N}$; the set $\Rep_{0}+j$ is called the fibre~$j$,
and the zero fibre is identified with $\Rep_{0}$. The distinction
is essential: a point of $E^{\N}_{u}(\alpha)$ lies in a bounded window
about $n\mapsto\alpha n$ at every time step $n\in\N$, whereas a
point of $E^{\N}(\alpha)$ need only return to the bounded window
for infinitely many time steps and may visit the reflecting boundary
between successive high excursions. Note that $E^{\N}(0)$ and $E^{\N}_{u}(0)$
are analogues of the radial and uniformly radial limit sets of a Kleinian
group, respectively; see \cite{BJ97}.

We call $E^{\N}\left(0\right)$ the \emph{recurrent} set. We also
have a disjoint union of $\alpha$-escaping sets, consisting of a
recurrent and a transient part: 
\[
E^{\N}(\alpha)=\left(E^{\N}(\alpha)\cap E^{\N}(0)\right)\,\dot{\cup}\,\left(E^{\N}(\alpha)\setminus E^{\N}(0)\right)
\]
The two central geometric quantities are intrinsic Poincaré exponents
of the $\N$-extension. The $\alpha$-escaping Poincaré exponent $\delta^{\N}(\alpha)$
is defined from words whose \emph{endpoint} at time $n$ lies in a
bounded window about $n\mapsto\alpha n$. We introduce in addition
the \emph{uniform $\alpha$-escaping Poincaré exponent} $\delta^{\N}_{u}(\alpha)$,
obtained by requiring the entire reflected path up to time $n$ to
remain in a bounded window about $k\mapsto\alpha k$. Their formal
definitions are given in \Cref{sec:basic-spectra}. The definition
of $\delta^{\N}(0)$ may be viewed as an analogue of the critical
exponent of the restricted Poincaré series for Kleinian groups introduced
by Falk and Stratmann in \cite{FS06}.

The uniform and transient spectra are governed by the uniform exponent:
for every $0<\alpha\le a_{+}$,
\begin{equation}
\dim_{H}E^{\N}_{u}(\alpha)=\dim_{H}\left(E^{\N}(\alpha)\setminus E^{\N}(0)\right)=\delta^{\N}_{u}(\alpha).\label{eq:intro-transient-uniform}
\end{equation}
For $\alpha=0$ we have 
\[
\dim_{H}E^{\N}_{u}(0)=\dim_{H}E^{\N}(0)=\delta^{\N}_{u}(0)=\delta^{\N}(0).
\]
These statements are proved in~\Cref{thm:uniform-N}. This result
may be viewed as an $\alpha$-escaping analogue of the Bishop--Jones
formula for the radial and uniformly radial limit sets of Kleinian
groups.

The i.o.\ $\alpha$-escaping set behaves differently. We show that
\begin{equation}
\dim_{H}E^{\N}(\alpha)=\delta^{\N}(\alpha),\qquad0\le\alpha\le a_{+}.\label{eq:intro-main-new}
\end{equation}
Moreover, for $0<\alpha<a_{+}$,
\begin{equation}
\dim_{H}\left(E^{\N}(\alpha)\cap E^{\N}(0)\right)=\min\left\{ \delta^{\N}(\alpha),\delta^{\N}(0)\right\} .\label{eq:intro-recurrent-dimension}
\end{equation}
See~\Cref{thm:exact-reflected-spectrum} and~\Cref{thm:recurrent-reflected-escape}.

Unlike the case $\alpha=0$, the inequality $\delta^{\N}_{u}(\alpha)\le\delta^{\N}(\alpha)$
may be strict for $\alpha>0$. We therefore define the \emph{strict-gain
set} by
\begin{align*}
\mathcal{A}^{\N} & \coloneqq\left\{ \alpha\in(0,a_{+}):\delta^{\N}(\alpha)>\delta^{\N}_{u}(\alpha)\right\} .
\end{align*}
It consists exactly of the positive escape rates at which i.o.\ escape
has larger dimension than uniform escape. Combining \cref{eq:intro-transient-uniform,eq:intro-main-new,eq:intro-recurrent-dimension}
gives, for $0<\alpha<a_{+}$,
\[
\delta^{\N}(\alpha)=\max\left\{ \min\left\{ \delta^{\N}(\alpha),\delta^{\N}(0)\right\} ,\delta^{\N}_{u}(\alpha)\right\} .
\]
Consequently, if $\alpha\in\mathcal{A}^{\N}$ then $\delta^{\N}(\alpha)\le\delta^{\N}(0)$.
In fact the inequality is then strict; see~\Cref{cor:strict-reflection-gain}.

The two dimension spectra also have sharply different behaviour at
zero. The i.o.\ spectrum is continuous on $[0,a_{+})$, whereas the
uniform spectrum is discontinuous at $0$ if and only if the maximal-dimension
drift $\alpha_{\max}$ defined in~\cref{eq:alpha-max} is negative.
Equivalently, this occurs if and only if there exists $\varepsilon>0$
such that
\begin{figure}[!ht]
\centering \begin{tikzpicture}
\begin{axis}[
    width=0.72\textwidth,
    height=0.23\textheight,
    xmin=0, xmax=1,
    ymin=0, ymax=1.05,
    axis lines=left,
        xtick={0,0.3912199,1},
    xticklabels={$0$,$\alpha^*$,$a_+=1$},
    ytick={0,0.5,1},
    tick label style={font=\small},
    label style={font=\small},
    legend style={draw=none, fill=none, font=\small, at={(0.96,0.96)}, anchor=north east},
]
\pgfmathdeclarefunction{HH}{1}{%
  \pgfmathparse{(#1<=0 || #1>=1) ? 0 : -(#1)*ln(#1) - (1-#1)*ln(1-#1)}%
}
\addplot[black,dashed,thick,domain=0:1,samples=220]
{HH((1+x)/2)/(-((1-x)/2)*ln(0.9)-((1+x)/2)*ln(0.1))};
\addlegendentry{$\dim_H E_{u}^{\N}(\alpha)$}
\addplot[black,solid,thick,domain=1:0.3761469014,samples=220,variable=\t]
({ln(exp(\t*ln(0.9))+exp(\t*ln(0.1)))/(\t*ln(9))},{\t});
\addlegendentry{$\dim_H E^{\N}(\alpha)$}
\addplot[black,solid,thick,domain=0.3912199288:1,samples=220]
{HH((1+x)/2)/(-((1-x)/2)*ln(0.9)-((1+x)/2)*ln(0.1))};
\addplot[black,dotted] coordinates {(0.3912199288,0) (0.3912199288,0.42)};
\end{axis}\end{tikzpicture} \caption{A locally constant example for which the strict-gain set is $\mathcal{A}^{\mathbb{N}}=(0,\alpha^{*})$,
where $\alpha^{*}\approx0.3912$. Thus the i.o.\ and uniform Hausdorff
dimension spectra differ exactly for $0<\alpha<\alpha^{*}$ and agree
at $\alpha=0$ and on $[\alpha^{*},1]$. The uniform spectrum is real
analytic on $(0,1)$ and discontinuous at $\alpha=0$. The i.o.\
spectrum is real analytic on $(0,1)\setminus\{\alpha^{*}\}$ and is
$C^{1}$ but not $C^{2}$ at $\alpha^{*}$.}
\label{fig:nonconcave-reflected}
\end{figure}
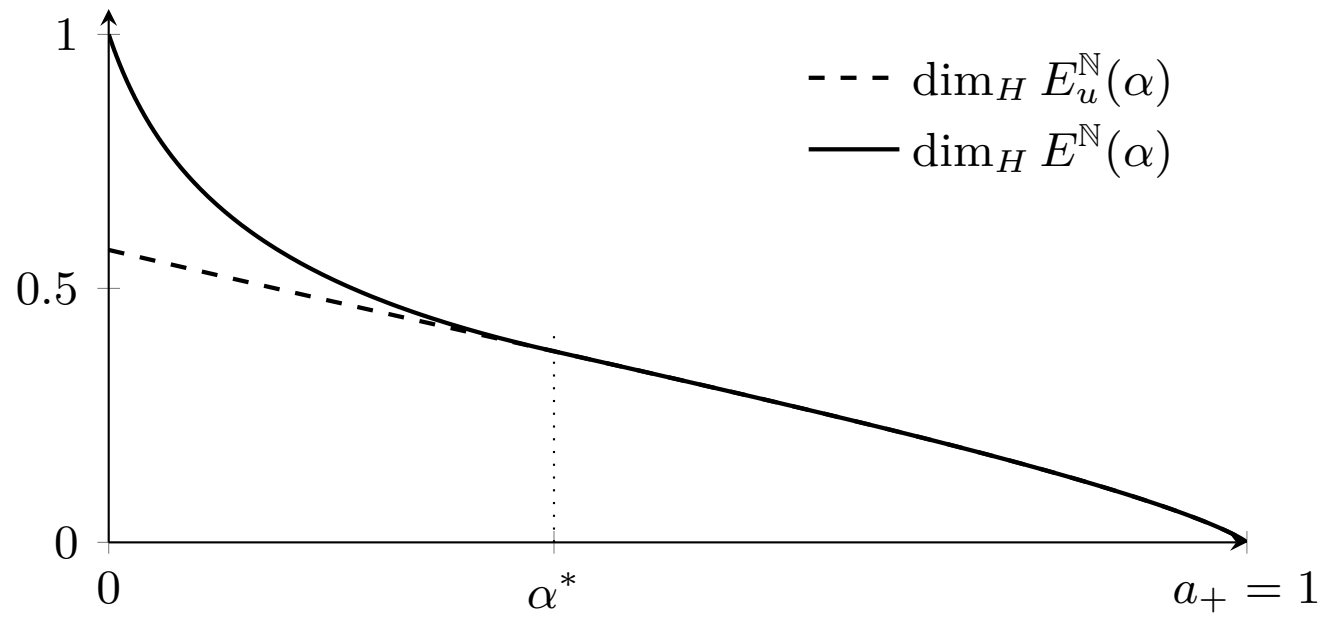
\[
(0,\varepsilon)\subset\mathcal{A}^{\N}.
\]
See~\Cref{cor:continuity-zero-jump}; an illustration, worked out
in~\Cref{ex:nonconcave-reflected}, is given in~\Cref{fig:nonconcave-reflected}.

A key computational ingredient is the \emph{$\alpha$-escaping pressure}
$\cP^{\N}_{\alpha}(g)$. Note that the $0$-escaping pressure coincides
with the Gurevich pressure of the lifted potential; see \Cref{rem:gurevich-zero-speed}.
For a Hölder potential $g$, let $\cP^{0}(g)$ be the base pressure
and define 
\[
G_{g}(q)\coloneqq\max\left\{ \cP^{0}(g+q\psi),\inf_{r\le0}\cP^{0}(g+r\psi)\right\} .
\]
Then~\Cref{thm:truncated-base-pressure} gives the following representation
as a concave Legendre transformation:
\[
\cP^{\N}_{\alpha}(g)=\widehat{G}_{g}(\alpha)\coloneqq\inf_{q\in\R}\left\{ G_{g}(q)-q\alpha\right\} .
\]
For the geometric potential $\phi=-\log|F'_{0}|$ write $p(s,q)\coloneqq\cP^{0}(s\phi+q\psi)$.
Then $\delta^{\N}(\alpha)$ is the zero of $s\mapsto\cP^{\N}_{\alpha}(s\phi)$.

The zero-speed formula may also be viewed as the exponential counterpart
of the classical drift trichotomy for reflected random walks. Indeed,
after normalizing $\cP^{0}(g)=0$ and writing $\mu_{g}$ for the equilibrium
state of $g$, convexity gives
\[
\begin{aligned}\cP^{\N}_{0}(g) & =\inf_{q\le0}\cP^{0}(g+q\psi),\\
\cP^{\N}_{0}(g)=0 & \quad\Longleftrightarrow\quad\int\psi\,\mathrm{d}\mu_{g}\le0,\\
\cP^{\N}_{0}(g)<0 & \quad\Longleftrightarrow\quad\int\psi\,\mathrm{d}\mu_{g}>0.
\end{aligned}
\]
Thus non-positive drift entails no exponential cost for returning
to the reflecting region, whereas positive drift makes such returns
exponentially rare. Lalley's local return theorem~\cite{Lalley95}
supplies, in the independent-increment setting, the corresponding
polynomial prefactors; the pressure records only the exponential rate.

For each $s$ let $q^{*}(s)$ be the nontrivial solution of $p(s,q)=p(s,0)$,
with $q^{*}(s)=0$ when $\partial_{q}p(s,0)=0$, and put $a^{*}(s):=\partial_{q}p(s,q^{*}(s))$.
The strict-gain criterion then takes the following intrinsic form
in terms of the base pressure alone (see~\Cref{cor:strict-reflection-gain}):
\[
\mathcal{A}^{\N}=\left\{ \alpha\in(0,a_{+}):\alpha<a^{*}(\delta^{\N}_{u}(\alpha))\right\} .
\]

A second main result gives a complete regularity classification of
$\alpha\mapsto\delta^{\N}(\alpha)$. The i.o.\ spectrum is $C^{1}$
on $(0,a_{+})$ and real analytic on each component of $\mathcal{A}^{\N}$
and on the interior of its complement. Moreover, $\mathcal{A}^{\N}$
is a finite union of intervals bounded away from the maximal escape
rate, that is, for some $\varepsilon>0$ we have $(a_{+}-\varepsilon,a_{+}]\cap\mathcal{A}^{\N}=\varnothing,$
and hence $\delta^{\N}=\delta^{\N}_{u}$ on a left neighbourhood of
$a_{+}$; see~\Cref{rem:no-gain-near-max-speed}. If $\alpha_{0}$
is a non-isolated boundary point of $\mathcal{A}^{\N}$ and the contact
function
\[
H(s)\coloneqq p(s,0)-q^{*}(s)a^{*}(s)
\]
has a zero of multiplicity $m$ at $s_{0}=\delta^{\N}(\alpha_{0})$,
then $m$ is odd and the i.o.\ spectrum is $C^{2m-1}$, but not $C^{2m}$,
at $\alpha_{0}$. In particular, a simple zero of $H$ produces a
second-order phase transition; see \Cref{thm:regularity-reflected-spectrum,cor:contact-multiplicity}.

\begin{figure}
\centering\begin{tikzpicture} \begin{groupplot}[ group style={group size=1 by 3, vertical sep=1.2cm}, width=0.78\textwidth, axis lines=left, xmin=0, xmax=17, xtick={0,2,4,6,8,10,12,14,16}, tick label style={font={\small }}, label style={font={\small }}, ] \nextgroupplot[ height=0.17\textheight, ymin=0, ymax=1.03,  legend pos=south west, legend style={draw=none,fill=none,font={\small }}, ] \addplot[black,dashed,thick,smooth] coordinates { (0,0.9924544) (0.5,0.9905506) (1,0.9883696) (1.5,0.9858891) (2,0.9830851) (2.5,0.9799310) (3,0.9763977) (3.5,0.9724535) (4,0.9680633) (4.5,0.9631887) (5,0.9577874) (5.5,0.9518127) (6,0.9452132) (6.5,0.9379321) (7,0.9299063) (7.3018,0.9246727) (8,0.9113352) (9,0.8888472) (10,0.8616367) (11,0.8287008) (12,0.7887621) (12.2437,0.7777860) (13,0.7400854) (14,0.6798702) (15,0.6015118) (15.5591,0.5417840) (16,0.4746313) (16.5,0.3415646) (16.7450,0.2246053) (16.8,0.1893176) (16.9,0.1114516) (16.9807,0.0274515) (17,0) }; \addlegendentry{$\dim_H E_{u}^{\N}(\alpha)$} \addplot[black,solid,thick,smooth] coordinates { (0,1.0000000) (0.5,0.9970253) (1,0.9938512) (1.5,0.9904580) (2,0.9868242) (2.5,0.9829251) (3,0.9787331) (3.5,0.9742171) (4,0.9693417) (4.5,0.9640673) (5,0.9583487) (5.5,0.9521352) (6,0.9453693) (6.5,0.9379862) (7,0.9299132) (7.3018,0.9246727) (8,0.9113352) (9,0.8888472) (10,0.8616367) (11,0.8287008) (12,0.7887621) (12.2437,0.7777860) (13,0.7401518) (14,0.6801959) (15,0.6017800) (15.5591,0.5417840) (16,0.4746313) (16.5,0.3415646) (16.7450,0.2246053) (16.8,0.1894700) (16.9,0.1121951) (16.9807,0.0274515) (17,0) }; \addlegendentry{$\dim_H E^{\N}(\alpha)$} \nextgroupplot[ height=0.17\textheight, ymin=0, ymax=0.0081,  scaled y ticks=false, ytick={0,0.002,0.004,0.006,0.008}, yticklabel style={/pgf/number format/fixed,/pgf/number format/precision=3}, ] \addplot[fill=gray!15,draw=none] coordinates {(16.32,0) (17,0) (17,0.0081) (16.32,0.0081)} --cycle; \addplot[black,thick,smooth] coordinates { (0,0.0075456) (0.5,0.0064748) (1,0.0054816) (1.5,0.0045689) (2,0.0037391) (2.5,0.0029941) (3,0.0023354) (3.5,0.0017636) (4,0.0012784) (4.5,0.0008786) (5,0.0005613) (5.5,0.0003225) (6,0.0001560) (6.5,0.0000541) (7,0.0000069) (7.3018,0) (12.2437,0) (13,0.0000664) (13.5,0.0001850) (14,0.0003257) (14.5,0.0003957) (15,0.0002683) (15.5,0.0000068) (15.5591,0) (16.7450,0) (16.7843,0.0000790) (16.8236,0.0002979) (16.8629,0.0005777) (16.9021,0.0007455) (16.9414,0.0005353) (16.9611,0.0002396) (16.9807,0) (17,0) }; \foreach\x in {7.3018,12.2437,15.5591,16.7450,16.9807}{ \addplot[black,dotted] coordinates {(\x,0) (\x,0.0081)}; } \nextgroupplot[ height=0.17\textheight, xmin=16.32, xmax=17, xtick={16.4,16.5,16.6,16.7,16.8,16.9,17}, ymin=0, ymax=0.00082,  scaled y ticks=false, ytick={0,0.0002,0.0004,0.0006,0.0008}, yticklabel style={/pgf/number format/fixed,/pgf/number format/precision=4}, ] \addplot[black,thick] coordinates { (16.32000,0.0000000) (16.34000,0.0000000) (16.36000,0.0000000) (16.38000,0.0000000) (16.40000,0.0000000) (16.42000,0.0000000) (16.44000,0.0000000) (16.46000,0.0000000) (16.48000,0.0000000) (16.50000,0.0000000) (16.52000,0.0000000) (16.54000,0.0000000) (16.56000,0.0000000) (16.58000,0.0000000) (16.60000,0.0000000) (16.62000,0.0000000) (16.64000,0.0000000) (16.66000,0.0000000) (16.68000,0.0000000) (16.70000,0.0000000) (16.72000,0.0000000) (16.74000,0.0000000) (16.74497,0.0000000) (16.74497,0.0000000) (16.74500,0.0000000) (16.74896,0.0000008) (16.75000,0.0000013) (16.75296,0.0000033) (16.75500,0.0000052) (16.75695,0.0000074) (16.76000,0.0000116) (16.76095,0.0000131) (16.76494,0.0000205) (16.76500,0.0000206) (16.76894,0.0000295) (16.77000,0.0000322) (16.77294,0.0000402) (16.77500,0.0000463) (16.77693,0.0000524) (16.78000,0.0000628) (16.78093,0.0000662) (16.78492,0.0000815) (16.78500,0.0000818) (16.78892,0.0000983) (16.79000,0.0001031) (16.79291,0.0001166) (16.79500,0.0001267) (16.79691,0.0001363) (16.80000,0.0001524) (16.80091,0.0001573) (16.80490,0.0001796) (16.80500,0.0001801) (16.80890,0.0002030) (16.81000,0.0002097) (16.81289,0.0002276) (16.81500,0.0002410) (16.81689,0.0002532) (16.82000,0.0002737) (16.82089,0.0002797) (16.82488,0.0003069) (16.82500,0.0003078) (16.82888,0.0003349) (16.83000,0.0003428) (16.83287,0.0003633) (16.83500,0.0003787) (16.83687,0.0003922) (16.84000,0.0004150) (16.84086,0.0004213) (16.84486,0.0004505) (16.84500,0.0004515) (16.84886,0.0004796) (16.85000,0.0004879) (16.85285,0.0005084) (16.85500,0.0005236) (16.85685,0.0005367) (16.86000,0.0005585) (16.86084,0.0005642) (16.86484,0.0005908) (16.86500,0.0005919) (16.86883,0.0006162) (16.87000,0.0006234) (16.87283,0.0006402) (16.87500,0.0006526) (16.87683,0.0006626) (16.88000,0.0006789) (16.88082,0.0006829) (16.88482,0.0007009) (16.88500,0.0007017) (16.88881,0.0007164) (16.89000,0.0007205) (16.89281,0.0007291) (16.89500,0.0007346) (16.89681,0.0007385) (16.90000,0.0007435) (16.90080,0.0007444) (16.90480,0.0007465) (16.90500,0.0007465) (16.90879,0.0007444) (16.91000,0.0007429) (16.91279,0.0007378) (16.91500,0.0007321) (16.91678,0.0007264) (16.92000,0.0007135) (16.92078,0.0007099) (16.92478,0.0006880) (16.92500,0.0006866) (16.92877,0.0006605) (16.93000,0.0006509) (16.93277,0.0006273) (16.93500,0.0006062) (16.93676,0.0005882) (16.94000,0.0005522) (16.94076,0.0005432) (16.94475,0.0004927) (16.94500,0.0004893) (16.94875,0.0004367) (16.95000,0.0004182) (16.95275,0.0003759) (16.95500,0.0003401) (16.95674,0.0003116) (16.96000,0.0002573) (16.96074,0.0002449) (16.96473,0.0001783) (16.96500,0.0001738) (16.96873,0.0001145) (16.97000,0.0000956) (16.97273,0.0000585) (16.97500,0.0000325) (16.97672,0.0000170) (16.98000,0.0000006) (17.00000,0.0000000) }; \foreach\x in {16.7450,16.9807}{ \addplot[black,dotted] coordinates {(\x,0) (\x,0.00082)}; } \end{groupplot} \end{tikzpicture}
\caption{A locally constant example for which the strict-gain set has three
components (see \Cref{ex:three-component-gain-set}). The upper panel
shows the i.o.\ and uniform Hausdorff dimension spectra on $[0,a_{+}]$,
the middle panel shows their difference, and the lower panel magnifies
this difference on the shaded interval. The uniform spectrum is real-analytic
on $(0,a_{+})$ and discontinuous at $0$. The i.o.\ spectrum is
real-analytic on $(0,a_{+})$ except at the five boundary points of
the components of the strict-gain set, where it is $C^{1}$ but not
$C^{2}$.}
\label{fig:three-component-gain-set}
\end{figure}
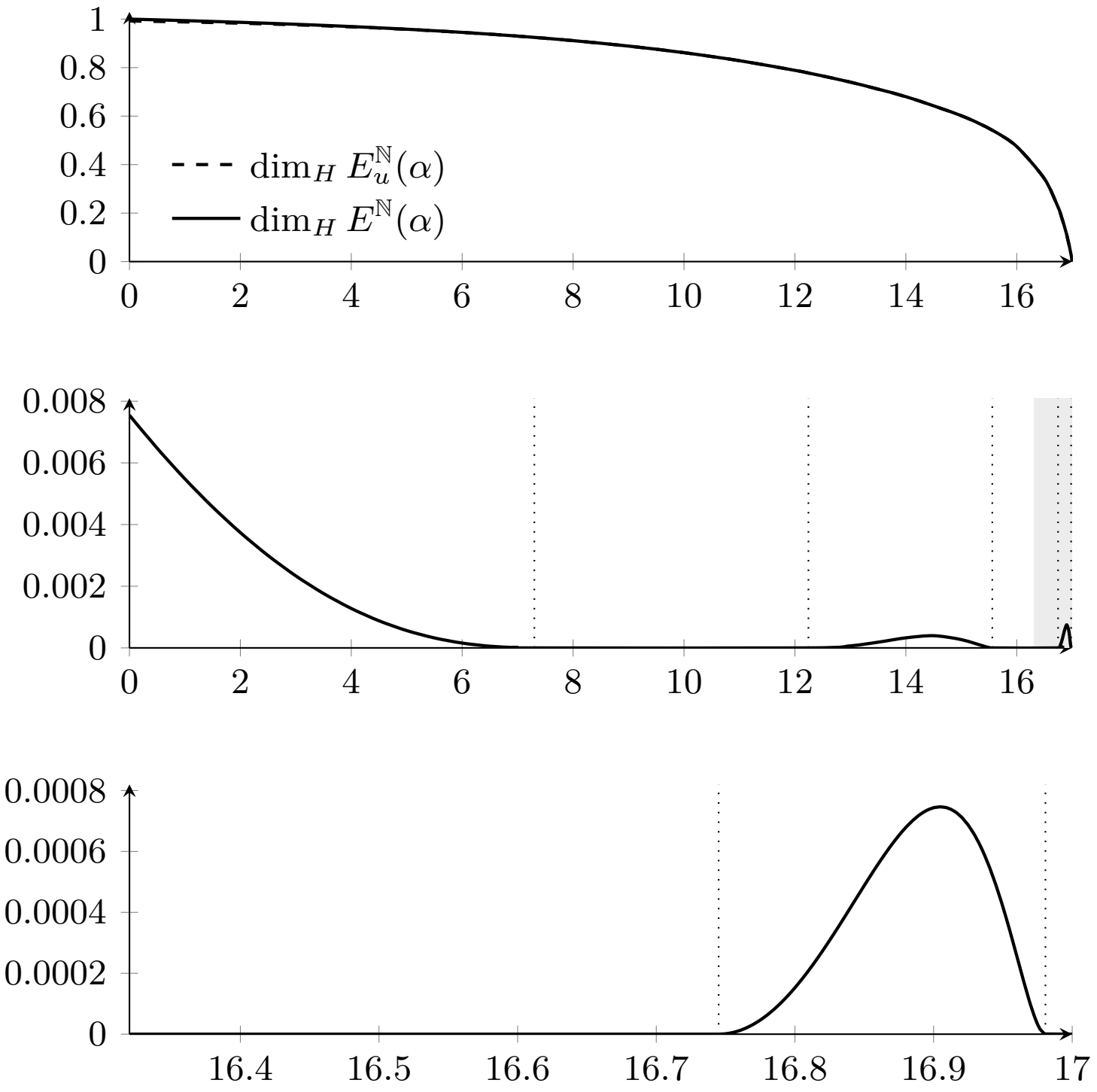

By contrast, $\alpha\mapsto\delta^{\N}_{u}(\alpha)$ is real analytic
on $(0,a_{+})$. The three-component example \Cref{ex:three-component-gain-set}
is illustrated in \Cref{fig:three-component-gain-set}. 

The lower bound for $E^{\N}(\alpha)\cap E^{\N}(0)$ is obtained from
a Moran construction built from a \emph{boundary waiting} phase, a
\emph{positive excursion}, and a \emph{downward recovery}. The necessary
uniform Frostman estimates for this construction are given in \Cref{sec:bounded-window-families}.

For comparison, the associated $\Z$-extension is defined by
\begin{equation}
F_{\Z}(x,k)=(F_{0}x,k+\psi(x)),\qquad(x,k)\in\Rep_{0}\times\Z.\label{eq:Z_extension}
\end{equation}
 Using as before the identification of $(x,k)$ with $x+k$, $F_{\Z}$
is the $\Z$-extension represented in \cite[Section~2.2]{GJK-reflected}
as a skew-periodic interval map without a reflecting boundary. By
\cite[Theorem 1.5]{GJK-realline} we have $\dim_{H}E^{\Z}(\alpha)=\dim_{H}E^{\Z}_{u}(\alpha)$
for every $\alpha$ in the attainable range (see \Cref{sec:basic-spectra}).
This is in stark contrast to the reflected $\N$-extension, for which
the two dimension spectra can differ. Thus, the separation between
i.o.\ and uniform escape is a genuine effect of the reflecting boundary.

The paper is organized as follows. In~\Cref{sec:setting} we set
up the symbolic model; in~\Cref{sec:basic-spectra} we define the
two Poincaré exponents of the $\N$-extension and record the basic
spectra; in~\Cref{sec:pressure} we introduce the $\alpha$-escaping
pressure representations, and discuss its connection to the Gurevich
pressure; in~\Cref{sec:bounded-window-families} we construct bounded-window
waiting, excursion and recovery families; in~\Cref{sec:general-reflected-pressure}
we prove the escape pressure, regularity and Hausdorff-dimension results;
and in~\Cref{sec:examples} we collect explicit locally constant
examples.

\section{Finite-branch skew products}

\label{sec:setting}

\subsection{The expanding base}

Recall the definition of the finite alphabet $I$ with $|I|\geq2$,
and $\Delta_{0}=\biguplus_{i\in I}J_{i}\subset(0,1)$ stated in the
introduction. We assume that $F_{0}:J_{i}\longrightarrow(0,1)$ is
a $C^{2}$ diffeomorphism for every $i\in I$, with uniform expansion
and bounded distortion. Its repeller is 
\[
\Rep_{0}\coloneqq\bigcap_{n\geq0}F^{-n}_{0}(\Delta_{0}).
\]
The symbolic model is the full shift $(\Sigma,\sigma)$ with $\Sigma=I^{\N}$
and coding map $\pi:\Sigma\to\Rep_{0}$. For a word $\omega=\omega_{0}\cdots\omega_{n-1}\in I^{n}$
we write $[\omega]$ for its cylinder and $J_{\omega}$ for the corresponding
basic interval.

The \emph{geometric potential }is 
\[
\phi\coloneqq-\log|F_{0}'|<0,
\]
which is Hölder. For $g:\Sigma\to\R$ write $S_{n}g=\sum^{n-1}_{j=0}g\circ\sigma^{j}$.

Let us emphasize that our assumptions force the step length function
$\psi:\Sigma\to\Z$ not to be cohomologous to a constant.

\subsection{Reflection at the boundary}

We use \emph{Lindley reflection} throughout:
\[
\rho_{i}(k)\coloneqq\max\{0,k+\psi(i)\},\qquad i\in I,\quad k\in\N.
\]
Thus, with $m_{-}\coloneqq-a_{-}$, the $\N$-extension agrees with
the $\Z$-extension $F_{\Z}$ as given in \cref{eq:Z_extension} on
every level~$k\geq m_{-}$. Recall that 
\[
F_{\N}(x,k)=(F_{0}x,\rho_{x_{0}}(k)).
\]
This is the $\N$-extension corresponding to the \emph{skew-periodic
interval map} with a left reflecting boundary in \cite[Section~2.2]{GJK-reflected},
here specialized to Lindley reflection.

For $\omega\in\Sigma$, define the height after $n$ steps starting
in fibre zero by
\begin{equation}
S^{\N}_{n}\psi(\omega):=h\bigl(F^{n}_{\N}(\pi(\omega))\bigr).\label{eq:NotationS^N}
\end{equation}
 Note that, despite the notation, $S^{\N}_{n}\psi$ is \textbf{not}
an additive cocycle because of the reflection at the boundary. Since
$\psi$ depends only on the first symbol, these heights depend only
on the first $n$ symbols of $\omega$. We therefore use the same
notation for finite words of length $n$.
\begin{lem}
\label{lem:transitive}$F_{\Z}$ and $F_{\N}$ are topologically transitive.
Moreover, for every $j\in\N$ there is a finite word $\tau$ whose
reflected path started at height $0$ ends at height $j$, that is,
$S^{\N}_{|\tau|}\psi(\tau)=j$.
\end{lem}

\begin{proof}
By Bézout's identity, every integer displacement can be written as
an integer combination of the allowed steps. Since steps of both signs
occur, one may add a positive zero-sum relation to make all coefficients
nonnegative. Hence every integer displacement can be realized by a
finite word, and therefore $F_{\Z}$ is topologically transitive.

For $F_{\N}$, every level can be brought to $0$ by repeating a negative
step. Conversely, for every $k\in\N$, choose a word with total displacement
$k$ and place all positive steps before the negative ones. Then the
corresponding path starting at $0$ never crosses the boundary and
ends at $k$. Hence $F_{\N}$ is topologically transitive, and the
word just constructed is the word $\tau$ of the last assertion.
\end{proof}

\section{Poincaré exponents and basic spectra}\label{sec:basic-spectra}

To unify the notation for the endpoint maps of the $\Z$- and $\N$-extensions,
we complement~\cref{eq:NotationS^N} by 
\[
S^{\Z}_{n}\psi(\omega)\coloneqq S_{n}\psi(\omega).
\]
 The attainable escape rate interval of the $\Z$-extension is 
\[
\cA=\left\{ \int\psi\,\dd\mu:\mu\in\cM_{\sigma}(\Sigma)\right\} =[a_{-},a_{+}],
\]
due to the fact that the base is the full shift and $\psi$ is constant
on one-cylinders. For $\mathbb{H}\in\{\Z,\N\}$, $K>0$, and $\alpha\in\mathcal{A}$
(with $\alpha\in\mathcal{A}_{+}\coloneqq\mathcal{A}\cap\R_{\geq0}$
if $\mathbb{H}=\N$), the \emph{$\alpha$-escaping Poincaré exponent
}is defined by 
\begin{equation}
{\delta^{\mathbb{H}}_{K}(\alpha)\coloneqq\inf\left\{ s\ge0:\sum_{n\ge1}\sum_{\substack{\omega\in I^{n}\\
|S^{\mathbb{H}}_{n}\psi(\omega)-\alpha n|\le K
}
}\exp\bigl(\sup_{[\omega]}S_{n}(s\phi)\bigr)<\infty\right\} .}\label{eq:deltaK-H}
\end{equation}
We then set 
\begin{equation}
{\delta^{\mathbb{H}}(\alpha)\coloneqq\lim_{K\to\infty}\delta^{\mathbb{H}}_{K}(\alpha),\qquad\mathbb{H}\in\{\Z,\N\}.}\label{eq:delta-H}
\end{equation}
For the $\N$-extension and $\alpha\in\mathcal{A}_{+}$, we also define
the \emph{uniform path family}
\[
\mathcal{U}^{\N}_{n,K}(\alpha)\coloneqq\left\{ \omega\in I^{n}:\max_{1\le k\le n}|S^{\N}_{k}\psi(\omega)-\alpha k|\le K\right\} .
\]
Its critical exponent is
\[
\delta^{\N}_{u,K}(\alpha)\coloneqq\inf\left\{ s\ge0:\sum_{n\ge1}\sum_{\omega\in\mathcal{U}^{\N}_{n,K}(\alpha)}\exp\bigl(\sup_{[\omega]}S_{n}(s\phi)\bigr)<\infty\right\} ,
\]
and we put
\[
\delta^{\N}_{u}(\alpha)\coloneqq\lim_{K\to\infty}\delta^{\N}_{u,K}(\alpha).
\]
Thus $\delta^{\N}(\alpha)$ only constrains the endpoint at time $n$,
whereas $\delta^{\N}_{u}(\alpha)$ constrains the entire reflected
path up to time $n$. We also write 
\[
\delta_{0}\coloneqq\dim_{H}\Rep_{0}.
\]

Recall the definition of $E^{\N}(\alpha)$ and $E^{\N}_{u}(\alpha)$
for $\alpha\geq0$. When replacing $\N$ by $\Z$, we define analogously
$E^{\Z}(\alpha)$ and $E^{\Z}_{u}(\alpha)$ for $\alpha\in\mathcal{A}$.

\emph{Reduction to the zero fibre.} All dimension statements below
may be verified on the zero fibre $\Rep_{0}$. Indeed, the fibre coordinates
of $x$ and $F_{\N}x$ differ by at most $\max\{a_{+},m_{-}\}$, so
the window conditions defining $E^{\N}(\alpha)$ and $E^{\N}_{u}(\alpha)$
hold for $x$ if and only if they hold for $F_{\N}x$, after enlarging
$K$; hence these sets, and with them $E^{\N}\left(0\right)$, $E^{\N}_{\mathrm{rec}}(\alpha)$
and $E^{\N}_{\mathrm{tr}}(\alpha)$, are completely $F_{\N}$-invariant,
that is, $F^{-1}_{\N}(E)=E$. For such a set $E$ the fibres $\Rep_{0}+j$,
$j\in\N$, cover $\Rep_{0}+\N$, so $\dim_{H}E=\sup_{j}\dim_{H}(E\cap(\Rep_{0}+j))$
by countable stability of Hausdorff dimension. By irreducibility (\Cref{lem:transitive})
there is a word $\tau\in I^{m}$ with $S^{\N}_{m}\psi(\tau)=j$; by
bounded distortion $F^{m}_{\N}$ maps $J_{\tau}\cap\Rep_{0}$ bi-Lipschitz
onto $\Rep_{0}+j$, and by complete invariance it maps $E\cap J_{\tau}$
onto $E\cap(\Rep_{0}+j)$. Therefore $\dim_{H}(E\cap(\Rep_{0}+j))\le\dim_{H}(E\cap\Rep_{0})$,
and consequently $\dim_{H}E=\dim_{H}(E\cap\Rep_{0})$; compare the
analogous arguments in~\cite{GJK-realline}.

The spectrum of the $\Z$-extension is known from \cite[Theorem~1.5]{GJK-realline}.
The cited theorem covers interior escape rates; at the two endpoints
the same identities follow by restricting to symbols with minimal
or maximal step.
\begin{thm}[The $\Z$-extension]
\label{thm:Z-spectrum} For every $\alpha\in\cA$, 
\[
{\dim_{H}E^{\Z}(\alpha)=\dim_{H}E^{\Z}_{u}(\alpha)=\delta^{\Z}(\alpha).}
\]
 Moreover, $\alpha\mapsto\delta^{\Z}(\alpha)$ is real analytic on
the interior of the attainable escape rate interval.
\end{thm}

\begin{prop}[Uniform Poincaré exponent and the associated $\Z$-extension]
\label{prop:uniform-exponent-comparison} We have
\[
\delta^{\N}_{u}(\alpha)=\begin{cases}
\delta^{\Z}(\alpha), & \alpha\in(0,a_{+}],\\
\delta^{\N}(0) & \alpha=0.
\end{cases}
\]
\end{prop}

\begin{proof}
Fix $\alpha>0$ and $K$. If a reflected path satisfies $|S^{\N}_{k}\psi-\alpha k|\le K$
for $1\le k\le n$, then it cannot meet the boundary after time $K/\alpha+1$.
Hence, by modifying only a bounded initial segment and enlarging the
window by a bounded amount, we can compare the uniform path families
of the $\N$-extension with those of the associated $\Z$-extension.
This proves $\delta^{\N}_{u}(\alpha)\le\delta^{\Z}(\alpha)$. The
reverse inequality follows in the same way. At $\alpha=0$, the condition
defining $\mathcal{U}^{\N}_{n,K}(0)$ confines the path to the finite
strip $\{0,\ldots,K\}$. The exhaustion principle for the Gurevich
pressure \cite[Lemma~4.5 and Proposition~2.1]{GJK-reflected} therefore
gives $\delta^{\N}_{u}(0)=\delta^{\N}(0)$.
\end{proof}

\begin{thm}[Uniform and transient spectra of the $\N$-extension]
\label{thm:uniform-N} We have
\[
\dim_{H}E^{\N}_{u}(\alpha)=\delta^{\N}_{u}(\alpha)=\begin{cases}
\dim_{H}\bigl(E^{\N}(\alpha)\setminus E^{\N}\left(0\right)\bigr), & \alpha\in(0,a_{+}],\\
\dim_{H}E^{\N}(0)=\delta^{\N}(0), & \alpha=0.
\end{cases}
\]
\end{thm}

\begin{proof}
Let $\alpha>0$. If $x\in E^{\N}(\alpha)\setminus E^{\N}(0)$, then
$h(F^{n}_{\N}x)\to\infty$, so reflection is inactive after finitely
many iterates. Hence, up to a finite prefix and a fibre translation,
the remaining orbit belongs to $E^{\Z}(\alpha)$. Therefore, by \Cref{thm:Z-spectrum}
and \Cref{prop:uniform-exponent-comparison}, 
\[
\dim_{H}\bigl(E^{\N}(\alpha)\setminus E^{\N}(0)\bigr)\le\dim_{H}\left(E^{\Z}(\alpha)\right)=\delta^{\N}_{u}(\alpha).
\]
Since $E^{\N}_{u}(\alpha)\subset E^{\N}(\alpha)\setminus E^{\N}(0)$,
it remains to show that $\dim_{H}E^{\N}_{u}(\alpha)\ge\delta^{\N}_{u}(\alpha)$.
Since $\delta^{\N}_{u}(\alpha)=\dim_{H}\left(E^{\Z}_{u}(\alpha)\right)$
by \Cref{thm:Z-spectrum} and \Cref{prop:uniform-exponent-comparison},
this follows by adding a finite prefix to the uniformly $\alpha$-escaping
orbits of $E^{\Z}_{u}(\alpha)$. This proves the first assertion.
For $\alpha=0$, the conclusion follows from \Cref{prop:uniform-exponent-comparison}
together with \cite[Proposition~2.1]{GJK-reflected}.
\end{proof}

\begin{thm}[Infinitely often spectrum of the $\N$-extension]
\label{thm:exact-reflected-spectrum} For every $\alpha\in\mathcal{A}_{+}$,
\begin{equation}
{\dim_{H}E^{\N}(\alpha)=\delta^{\N}(\alpha).}\label{eq:exact-reflected-spectrum}
\end{equation}
For $0<\alpha<a_{+}$, 
\begin{equation}
{\dim_{H}\bigl(E^{\N}(\alpha)\cap E^{\N}(0)\bigr)=\min\{\delta^{\N}(\alpha),\delta^{\N}(0)\}.}\label{eq:main-recurrent-dimension}
\end{equation}
The component realizing the full dimension can be identified from
$\delta^{\N}_{u}(\alpha)$ and $\delta^{\N}(0)$ alone. More precisely,
\begin{align}
\delta^{\N}_{u}(\alpha)\le\delta^{\N}(0) & \quad\Longrightarrow\quad\dim_{H}\bigl(E^{\N}(\alpha)\cap E^{\N}(0)\bigr)=\delta^{\N}(\alpha),\label{eq:recurrent-realizes}\\[1mm]
\delta^{\N}_{u}(\alpha)>\delta^{\N}(0) & \quad\Longrightarrow\quad\dim_{H}E^{\N}(\alpha)\setminus E^{\N}\left(0\right)=\delta^{\N}_{u}(\alpha)=\delta^{\N}(\alpha).\label{eq:translation-realizes}
\end{align}
\end{thm}

The proof of \Cref{thm:exact-reflected-spectrum} is given in \Cref{sec:general-reflected-pressure},
after the pressure representations and the bounded-window estimates
needed for the lower bounds.


\section{Pressure representations of the Poincaré exponents}

\label{sec:pressure} 

We now introduce exponential growth rates as tools for evaluating
the Poincaré exponents defined in \eqref{eq:delta-H}. The geometric
quantities $\delta^{\mathbb{H}}(\alpha)$ remain the primary definitions
throughout.

\subsection{Escape pressures and the return case}

For $\mathbb{H}\in\{\Z,\N\}$, $K>0$, a Hölder potential $g$, and
$\alpha$, put 
\[
Z^{\mathbb{H}}_{n,K}(g,\alpha)\coloneqq\sum_{\substack{\omega\in I^{n}\\
|S^{\mathbb{H}}_{n}\psi(\omega)-\alpha n|\le K
}
}\exp\!\left(\sup_{\tau\in[\omega]}S_{n}g(\tau)\right)
\]
and define 
\[
\cP^{\mathbb{H}}_{\alpha}(g)\coloneqq\lim_{K\to\infty}\limsup_{n\to\infty}\frac{1}{n}\log Z^{\mathbb{H}}_{n,K}(g,\alpha),\qquad\mathbb{H}\in\{\Z,\N\}.
\]
The base pressure has no escape rate parameter and is denoted by
\[
\cP^{0}(g)\coloneqq\lim_{n\to\infty}\frac{1}{n}\log\sum_{\omega\in I^{n}}\exp\bigl(\sup_{[\omega]}S_{n}g\bigr).
\]
 For the geometric potential we abbreviate 
\[
p^{\mathbb{H}}_{\alpha}(s)\coloneqq\cP^{\mathbb{H}}_{\alpha}(s\phi).
\]
 Thus $p^{\Z}_{v}(s)$ is the exponential rate of words of the $\Z$-extension
with endpoint $vn+O(1)$, while $p^{\N}_{\alpha}(s)$ is the corresponding
rate for words of the $\N$-extension ending at height $\alpha n+O(1)$.
For fixed $K$, the inner limsup is again a $1$-induced pressure
in the sense of \cite[Definition~1.1]{JKL14}, now for the word collection
cut out by the corresponding bounded escape rate window.
\begin{prop}[Root representation of the Poincaré exponents]
\label{prop:delta-pressure-root} For every $\alpha\in\mathcal{A}$
if $\mathbb{H}=\Z$ and every $\alpha\in\mathcal{A}_{+}$ if $\mathbb{H}=\N$,
\begin{equation}
{\delta^{\mathbb{H}}(\alpha)=\sup\{s\ge0:p^{\mathbb{H}}_{\alpha}(s)\ge0\},\qquad\mathbb{H}\in\{\Z,\N\}.}\label{eq:deltaH-marked-root}
\end{equation}
\end{prop}

\begin{proof}
For fixed $K$ the root test identifies the abscissa of convergence
in \eqref{eq:deltaK-H} with the zero threshold of the corresponding
limsup exponential rate. This is the fixed-length specialization of
the critical-exponent principle in \cite[Theorem~2.4]{JKL14}. The
bounded-window exponential rates are non-decreasing in $K$ and converge
to $p^{\mathbb{H}}_{\alpha}(s)$. Moreover, uniform negativity of
$\phi$ makes every such rate strictly decreasing in $s$ with a common
slope bound. Hence the zero thresholds converge to the zero threshold
of the limiting pressure. Letting $K\to\infty$ gives \eqref{eq:deltaH-marked-root}.
\end{proof}

Since $\phi$ is uniformly negative, there is $\kappa>0$ with $\phi\le-\kappa$.
Hence, whenever the pressures are finite and $t>s$, $p^{\mathbb{H}}_{\alpha}(t)\le p^{\mathbb{H}}_{\alpha}(s)-\kappa(t-s)$,
so $s\mapsto p^{\mathbb{H}}_{\alpha}(s)$ is strictly decreasing.
In particular, $s<\delta^{\mathbb{H}}(\alpha)$ implies $p^{\mathbb{H}}_{\alpha}(s)>0$.

The $\alpha$-escape pressure of the $\Z$-extension has the standard
constrained-pressure representation.
\begin{prop}[Escape pressure of the $\Z$-extension]
\label{prop:translation-marked-pressure} For every Hölder potential
$g$ and every $v\in\cA$, 
\begin{equation}
{\cP^{\Z}_{v}(g)=\inf_{q\in\R}\bigl\{\cP^{0}(g+q\psi)-qv\bigr\}.}\label{eq:pZ-pressure-formula-general}
\end{equation}
In particular, 
\[
p^{\Z}_{v}(s)=\inf_{q\in\R}\cP^{0}\bigl(s\phi+q(\psi-v)\bigr).
\]
Consequently, 
\[
\delta^{\Z}(v)=\inf_{q\in\R}\inf\left\{ s\geq0:\cP^{0}\bigl(s\phi+q(\psi-v)\bigr)\leq0\right\} .
\]
\end{prop}

\begin{proof}
For interior escape rates this is the finite-alphabet specialization
of \cite[Theorem~4.1]{GJK-realline}. The endpoint formulas follow
by the finite-alphabet limiting argument (equivalently, by restricting
to symbols with extremal step). The last formula follows from \Cref{prop:delta-pressure-root}.
\end{proof}

\begin{prop}[The $0$-escaping pressure and exact returns]
\label{prop:zero-speed-pressure} For a Hölder potential $g$ and
$\mathbb{H}\in\{\Z,\N\}$, put
\[
Z^{\mathbb{H}}_{n,0}(g)\coloneqq\sum_{\omega\in\mathcal{R}^{\mathbb{H}}_{n,0}}\exp\!\left(\sup_{\tau\in[\omega]}S_{n}g(\tau)\right)\;\text{with }\;\mathcal{R}^{\mathbb{H}}_{n,0}\coloneqq\left\{ \omega\in I^{n}:S^{\mathbb{H}}_{n}\psi(\omega)=0\right\} .
\]
Then the $0$-escaping pressure has the exact-return representation
\begin{equation}
{\cP^{\mathbb{H}}_{0}(g)=\limsup_{n\to\infty}\frac{1}{n}\log Z^{\mathbb{H}}_{n,0}(g).}\label{eq:zero-speed-exact-return}
\end{equation}
Moreover,
\begin{equation}
{\cP^{\Z}_{0}(g)=\inf_{q\in\R}\cP^{0}(g+q\psi),\qquad\cP^{\N}_{0}(g)=\inf_{q\le0}\cP^{0}(g+q\psi).}\label{eq:zero-speed-pressure-formulas}
\end{equation}
Consequently, $\delta^{\mathbb{H}}(0)=\sup\{s\ge0:p^{\mathbb{H}}_{0}(s)\ge0\}.$
\end{prop}

\begin{proof}
For fixed $K$, the condition $|S^{\mathbb{H}}_{n}\psi|\le K$ restricts
the endpoint to a finite set independent of $n$. Transitivity supplies
connectors of uniformly bounded length from each of these states to
the reference state $0$ and conversely. Bounded distortion changes
the corresponding partition sums only by multiplicative constants
and bounded shifts of length. Therefore the bounded windows at escape
rate $0$ and the exact-return families have the same exponential
rate, which proves \eqref{eq:zero-speed-exact-return}. The variational
formulas for the exact-return rates are \cite[Theorem~4.6(2)]{GJK-reflected}
for the $\Z$-extension and \cite[Theorem~4.13]{GJK-reflected} for
the $\N$-extension. This gives \eqref{eq:zero-speed-pressure-formulas}.
The last assertion is the case $\alpha=0$ of \Cref{prop:delta-pressure-root}.
\end{proof}

\begin{rem}
The earlier work \cite{GJK-reflected} considered only returns to
a reference fibre, so no escape-rate index was needed there. In the
present notation, those pressures are $\cP^{\mathbb{H}}_{0}(g)$.
The corresponding exponent is $\delta^{\mathbb{H}}(0)$, and we introduce
no separate return-pressure or return-exponent notation.
\end{rem}

Since $\delta_{0}=\dim_{H}\Rep_{0}$, Bowen's formula gives 
\[
\cP^{0}(\delta_{0}\phi)=0.
\]
Let $\mu_{\delta_{0}\phi}$ denote the equilibrium state of $\delta_{0}\phi$.
We define the maximal-dimension drift 
\begin{equation}
\alpha_{\max}\coloneqq\int\psi\,\dd\mu_{\delta_{0}\phi}.\label{eq:alpha-max}
\end{equation}
The terminology reflects the fact that $\delta^{\Z}(\alpha)$ reaches
its maximum $\delta_{0}$ at $\alpha=\alpha_{\max}$; see \cite[Proposition~4.3]{GJK-realline}.

\begin{rem}[Gurevich pressure, finite irreducibility, and compact exhaustion]
\label{rem:gurevich-zero-speed} For $\mathbb{H}\in\{\Z,\N\}$, let
$(\widehat{\Sigma}_{\mathbb{H}},\widehat{\sigma}_{\mathbb{H}})$ denote
the natural countable Markov coding of the $\mathbb{H}$-extension,
with alphabet $I\times\mathbb{H}$, and let $\widetilde{g}$ be the
lift of $g$ obtained by forgetting the fibre coordinates. The exact-return
representation in \Cref{prop:zero-speed-pressure} identifies $\cP^{\mathbb{H}}_{0}(g)$
with the Gurevich pressure of $\widetilde{g}$ in Sarig's thermodynamic
formalism \cite{Sarig1999}. Indeed, exact-return words give periodic
words of the extension, and bounded distortion and irreducibility
identify their exponential growth rate with that of periodic words
based at any fixed symbol. A possible period is accommodated by taking
a limsup.

The extensions are irreducible, but they are \emph{not finitely irreducible}
in the sense of Mauldin--Urbański~\cite[Section~2.1]{MU2003}: each
transition changes the fibre coordinate by at most $\max_{i\in I}|\psi(i)|$,
so no finite family of connecting words can join arbitrarily distant
fibres, and the connectors used in the proof of~\Cref{prop:zero-speed-pressure}
only serve a fixed finite set of endpoint heights. Consequently, the
comparison between the Gurevich pressure and the unrestricted cylinder-sum
pressure that is available under finite irreducibility~\cite[Section~7]{MU2001}
cannot be invoked here; in fact, the unrestricted cylinder sums of
$\widetilde{g}$ are infinite, because every base word has a lift
from every initial fibre. Instead, the exhausting principle of~\cite{JKL14}
identifies $\cP^{\mathbb{H}}_{0}(g)$ with the supremum of the classical
pressures on compact invariant subsystems: if $\widehat{\Sigma}_{\mathbb{H},L}$
denotes the finite-state subsystem of trajectories whose fibre coordinates
remain in $\mathbb{H}\cap[-L,L]$, then
\[
\cP^{\mathbb{H}}_{0}(g)=\sup_{L\ge1}P_{\mathrm{top}}\!\left(\widehat{\sigma}_{\mathbb{H}}\big|_{\widehat{\Sigma}_{\mathbb{H},L}},\widetilde{g}\big|_{\widehat{\Sigma}_{\mathbb{H},L}}\right),
\]
where $P_{\mathrm{top}}$ is the usual pressure on a compact subsystem.
\end{rem}

\section{Bounded-window word families}

\label{sec:bounded-window-families}

The lower bounds below require families whose centred displacement
stays in a fixed window at every intermediate time. We obtain them
directly from the Moran construction in~\cite[proof of Theorem~1.5]{GJK-realline}.
The only additional observation needed here is that the construction
uses one fixed main block length and a finite set of bounded correction
words. This gives uniform control between construction times and,
after normalizing the potential, a prefix Frostman estimate with an
arbitrarily small exponential loss. No local limit theorem, endpoint
conditioning or cycle lemma is used.

\subsection{Bounded-window Moran blocks}
\begin{lem}[Bounded-window Frostman blocks]
\label{lem:bounded-window-blocks} Let $g$ be Hölder, $v\in(a_{-},a_{+})$
and $\eta>0$. Then there are constants $K,C<\infty$, a compact set
$X=X(g,v,\eta)\subset\Sigma$ and a probability measure $\lambda$
supported on $X$ such that every $\omega\in X$ satisfies
\begin{equation}
|S_{n}\psi(\omega)-vn|\le K\qquad(n\ge0),\label{eq:bounded-window-tracking}
\end{equation}
and, for every $n\ge1$,
\begin{equation}
\lambda([\omega|_{n}])\le C\exp\left\{ \sup_{x\in[\omega|_{n}]}S_{n}g(x)-n\bigl(\cP^{\Z}_{v}(g)-\eta\bigr)\right\} .\label{eq:bounded-window-frostman}
\end{equation}
\end{lem}

\begin{proof}
Set $\psi_{v}\coloneqq\psi-v$, $a\coloneqq\cP^{\Z}_{v}(g)-\eta$,
and $\widetilde{g}\coloneqq g-a$. Then $\cP^{\Z}_{v}(\widetilde{g})=\eta>0$.
Since $\psi_{v}$ is constant on one-cylinders and $0\in(\min\psi_{v},\max\psi_{v})$,
the fibre-induced pressure in~\cite[Theorem~3.1 and Remark~3.7]{GJK-realline}
agrees, for every fixed tolerance $\varepsilon>0$, with the $v$-escape
pressure $\cP^{\Z}_{v}$. Hence
\[
\limsup_{n\to\infty}\frac{1}{n}\log\sum_{\substack{\nu\in I^{n}:\ |S_{n}\psi_{v}(\nu)|\le\varepsilon}
}\exp\bigl(\sup_{[\nu]}S_{n}\widetilde{g}\bigr)=\eta.
\]
In particular, for every $M>0$ there is a single sufficiently large
length $\ell$ for which, with
\[
\Gamma\coloneqq\{\nu\in I^{\ell}:|S_{\ell}\psi_{v}(\nu)|\le\varepsilon\},
\]
we have
\begin{equation}
Z_{\Gamma}\coloneqq\sum_{\nu\in\Gamma}\exp\bigl(\sup_{[\nu]}S_{\ell}\widetilde{g}\bigr)>M.\label{eq:Gamma-large-weight}
\end{equation}
This is precisely the starting point of the Moran construction in
the proof of Theorem~1.5 of the cited paper.

Apply \cite[Lemma~3.8]{GJK-realline} with $N=2\varepsilon$. It gives
a finite collection $\Lambda$ of correction words and $r<\infty$
such that every word whose centred displacement has absolute value
at most $2\varepsilon$ can be followed by some $\tau\in\Lambda$,
with $|\tau|\le r$, so that the corrected displacement again has
absolute value at most $\varepsilon$. Choose one such correction
for each parent--main-block pair. Starting with words in $\Gamma$
and repeating the same family $\Gamma$ at every generation, followed
by such a correction, gives a rooted tree of words. At every tree
level
\[
|S_{m}\psi_{v}|\le\varepsilon,
\]
and two successive tree levels differ in length by a number in $[\ell,\ell+r]$.
Thus, with $L\coloneqq\ell+r$,
\[
|S_{n}\psi_{v}|\le\varepsilon+L\|\psi_{v}\|_{\infty}
\]
at every intermediate time. This proves \eqref{eq:bounded-window-tracking}
after increasing $K$.

We next put a measure on this tree. For every tree word $u$ set
\[
W(u)\coloneqq\exp\bigl(\sup_{[u]}S_{|u|}\widetilde{g}\bigr).
\]
The Bowen property of $\widetilde{g}$ and the finiteness of $\Lambda$
give a constant $A\ge1$, independent of the tree level, the parent
$u$ and the chosen main-block length $\ell$, such that
\[
\sum_{w\in\operatorname{Ch}(u)}W(w)\ge A^{-1}Z_{\Gamma}W(u).
\]
Indeed, a child $\omega=u\nu\eta$ of $u$ in $\operatorname{Ch}(u)$
is obtained by appending one $\nu\in\Gamma$ and one word $\eta$
from the fixed finite correction family; bounded distortion and the
bounded correction length change the product weight only by a uniform
factor. Choose the length $\ell$ in \eqref{eq:Gamma-large-weight}
so that $Z_{\Gamma}\ge A$. Define the transition probabilities by
\[
\lambda([w]\mid[u])\coloneqq\frac{W(w)}{\sum_{z\in\operatorname{Ch}(u)}W(z)}.
\]
Then
\begin{equation}
\lambda([w]\mid[u])\le\frac{W(w)}{W(u)}.\label{eq:tree-transition-frostman}
\end{equation}
The nested tree cylinders define a compact set $X$, and these transition
probabilities define a Borel probability measure $\lambda$ on $X$.
If $w$ is a tree word, multiplying \eqref{eq:tree-transition-frostman}
along the branch from the root to $w$ makes the intermediate weights
telescope and gives
\[
\lambda([w])\le W(w)=\exp\bigl(\sup_{[w]}S_{|w|}\widetilde{g}\bigr).
\]
Let $n\ge1$ be arbitrary and let $t^{-}(n)$ be the last tree time
not larger than $n$, with the root time $0$ allowed. Since consecutive
tree times differ by at most $L$,
\[
\lambda([\omega|_{n}])\le\lambda([\omega|_{t^{-}(n)}])\le\exp\bigl(\sup_{[\omega|_{t^{-}(n)}]}S_{t^{-}(n)}\widetilde{g}\bigr).
\]
The terminal piece has length at most $L$. The Bowen property of
$\widetilde{g}$ and boundedness of $\widetilde{g}$ therefore give
a constant $C_{0}$, independent of $n$, such that
\[
\lambda([\omega|_{n}])\le C_{0}\exp\bigl(\sup_{[\omega|_{n}]}S_{n}\widetilde{g}\bigr).
\]
Since $\widetilde{g}=g-a$ and $a=\cP^{\Z}_{v}(g)-\eta$, this proves
\eqref{eq:bounded-window-frostman}.
\end{proof}

Summing \eqref{eq:bounded-window-frostman} over the supported $n$-cylinders
and using that their $\lambda$-masses sum to one gives, for every
$n\ge1$,
\[
\sum_{\substack{\omega\in I^{n}:\\[1mm]{}
[\omega]\cap X\ne\varnothing
}
}\exp\!\left(\sup_{[\omega]}S_{n}g\right)\ge C^{-1}\exp\left(n(\cP^{\Z}_{v}(g)-\eta)\right).
\]

\subsection{Waiting, excursion and recovery families}

For a Hölder potential $g$ let $q^{\N}_{g}\le0$ be the minimizer
of
\[
q\longmapsto\cP^{0}(g+q\psi)\quad\text{on }(-\infty,0],
\]
and put $v^{\N}_{g}\coloneqq\left.\partial_{q}\cP^{0}(g+q\psi)\right|_{q=q^{\N}_{g}}.$
Then
\begin{equation}
v^{\N}_{g}\le0,\qquad q^{\N}_{g}v^{\N}_{g}=0,\qquad\cP^{\N}_{0}(g)=\cP^{\Z}_{v^{\N}_{g}}(g).\label{eq:complementarity-waiting}
\end{equation}
Indeed, if $q^{\N}_{g}<0$ the derivative vanishes, while if $q^{\N}_{g}=0$
the constrained minimum implies $v^{\N}_{g}\le0$. In either case
$q^{\N}_{g}$ is the minimizing parameter in the $\Z$-extension formula~\eqref{eq:pZ-pressure-formula-general}
at $v=v^{\N}_{g}$. Since the equilibrium state of $g+q\psi$ has
full support for every $q\in\R$, $v^{\N}_{g}\in(a_{-},a_{+})$. For
$g=s\phi$ we write $q^{\N}_{s}\coloneqq q^{\N}_{s\phi}$ and $v^{\N}_{s}\coloneqq v^{\N}_{s\phi}$.
\begin{prop}[Bounded-window phase families]
\label{prop:bounded-window-phase-families} Fix a Hölder potential
$g$, $\eta>0$, $\beta\in(0,a_{+})$ and $\gamma\in(0,-a_{-})$.
For $v\in\{v^{\N}_{g},\beta,-\gamma\}$ let $(X_{v},\lambda_{v},K_{v},C_{v})$
be the objects supplied by \Cref{lem:bounded-window-blocks} for $(g,v,\eta)$,
and define
\begin{equation}
\mathcal{B}_{n}(v)\coloneqq\{\omega\in I^{n}:\lambda_{v}([\omega])>0\},\qquad\lambda_{v,n}(\omega)\coloneqq\lambda_{v}([\omega]).\label{eq:Bv-definition}
\end{equation}
Choose $i_{+}\in I$ with $\psi(i_{+})=a_{+}$, choose $r_{+}\ge1$
so that $H_{+}\coloneqq r_{+}a_{+}>K_{\beta}$, and put $\tau^{+}\coloneqq i^{r_{+}}_{+}$.
For $n\ge r_{+}$ set
\[
\begin{array}{rlcrl}
\mathcal{W}_{n}(g) & \coloneqq\mathcal{B}_{n}(v^{\N}_{g}), &  & \lambda^{{\rm W}}_{n}(\omega) & \coloneqq\lambda_{v^{\N}_{g},n}(\omega),\\
\mathcal{U}_{n}(g,\beta) & \coloneqq\{\tau^{+}\omega:\omega\in\mathcal{B}_{n-r_{+}}(\beta)\}, &  & \lambda^{{\rm U}}_{n}(\tau^{+}\omega) & \coloneqq\lambda_{\beta,n-r_{+}}(\omega),\\
\mathcal{D}_{n}(g,\gamma) & \coloneqq\mathcal{B}_{n}(-\gamma), &  & \lambda^{{\rm D}}_{n}(\omega) & \coloneqq\lambda_{-\gamma,n}(\omega).
\end{array}
\]
For a prefix $\rho$ of length at most $n$, $\lambda^{\star}_{n}([\rho])$
denotes the total $\lambda^{\star}_{n}$-mass of words in the corresponding
family which begin with $\rho$. When $g=s\phi$ we abbreviate $\mathcal{W}_{n}(s)\coloneqq\mathcal{W}_{n}(s\phi)$,
$\mathcal{U}_{n}(s,\beta)\coloneqq\mathcal{U}_{n}(s\phi,\beta)$ and
$\mathcal{D}_{n}(s,\gamma)\coloneqq\mathcal{D}_{n}(s\phi,\gamma)$.
For a word $\omega$ and an initial height $h\ge0$ write
\[
h^{\N}_{k}(\omega;h)\coloneqq h+S_{k}\psi(\omega)-\min\!\left\{ 0,\min_{0\le m\le k}(h+S_{m}\psi(\omega))\right\} .
\]
\begin{enumerate}
\item There are constants $K_{{\rm W}},K_{{\rm U}},K_{{\rm D}}<\infty$
such that the following hold for every word in the indicated family
and every $0\le k\le n$. For $\omega\in\mathcal{W}_{n}(g)$,
\[
|S_{k}\psi(\omega)-v^{\N}_{g}k|\le K_{{\rm W}},\qquad0\le h^{\N}_{k}(\omega;h)\le h+2K_{{\rm W}}\quad(h\ge0).
\]
For $\omega\in\mathcal{U}_{n}(g,\beta)$,
\[
S_{k}\psi(\omega)\ge0,\qquad|S_{k}\psi(\omega)-\beta k|\le K_{{\rm U}}.
\]
For $\omega\in\mathcal{D}_{n}(g,\gamma)$,
\[
|S_{k}\psi(\omega)+\gamma k|\le K_{{\rm D}}.
\]
Consequently, if a recovery word is started at height $H\ge\gamma n+K_{{\rm D}}$,
reflection is inactive throughout it and
\[
\bigl|h^{\N}_{n}(\omega;H)-(H-\gamma n)\bigr|\le K_{{\rm D}}.
\]
\item There is $C<\infty$ such that, for every supported word and every
$1\le k\le n$,
\begin{align}
\lambda^{{\rm W}}_{n}([\omega|_{k}]) & \le C\exp\!\left(\sup_{[\omega|_{k}]}S_{k}g-k\bigl(\cP^{\N}_{0}(g)-\eta\bigr)\right),\label{eq:waiting-prefix-law}\\
\lambda^{{\rm U}}_{n}([\omega|_{k}]) & \le C\exp\!\left(\sup_{[\omega|_{k}]}S_{k}g-k\bigl(\cP^{\Z}_{\beta}(g)-\eta\bigr)\right),\label{eq:excursion-prefix-law}\\
\lambda^{{\rm D}}_{n}([\omega|_{k}]) & \le C\exp\!\left(\sup_{[\omega|_{k}]}S_{k}g-k\bigl(\cP^{\Z}_{-\gamma}(g)-\eta\bigr)\right).\label{eq:recovery-prefix-law}
\end{align}
\item After increasing $C$, for all sufficiently large $n$,
\begin{align*}
\sum_{\omega\in\mathcal{W}_{n}(g)}\e^{\sup_{[\omega]}S_{n}g} & \ge C^{-1}\e^{n(\cP^{\N}_{0}(g)-\eta)},\\
\sum_{\omega\in\mathcal{U}_{n}(g,\beta)}\e^{\sup_{[\omega]}S_{n}g} & \ge C^{-1}\e^{n(\cP^{\Z}_{\beta}(g)-\eta)},\\
\sum_{\omega\in\mathcal{D}_{n}(g,\gamma)}\e^{\sup_{[\omega]}S_{n}g} & \ge C^{-1}\e^{n(\cP^{\Z}_{-\gamma}(g)-\eta)}.
\end{align*}
\end{enumerate}
\end{prop}

\begin{proof}
For $\mathcal{W}_{n}(g)$ and $\mathcal{D}_{n}(g,\gamma)$, all tracking
and measure estimates are the corresponding assertions of \Cref{lem:bounded-window-blocks},
restricted to length-$n$ prefixes. The finite laws in \eqref{eq:Bv-definition}
are probability laws because $\lambda_{v}$ is supported on $X_{v}$.
For waiting, $v^{\N}_{g}\le0$ and
\[
S_{k}\psi-S_{m}\psi=v^{\N}_{g}(k-m)+e_{k}-e_{m}\le2K_{v^{\N}_{g}}\qquad(m\le k),
\]
so the Lindley formula gives $h^{\N}_{k}(\omega;h)\le\max\{h+K_{v^{\N}_{g}},2K_{v^{\N}_{g}}\}\le h+2K_{v^{\N}_{g}}$.
Thus one may take $K_{{\rm W}}=K_{v^{\N}_{g}}$. Identity \eqref{eq:complementarity-waiting}
turns the $\Z$-extension rate into $\cP^{\N}_{0}(g)$.

For excursions, if $\omega=\tau^{+}\nu\in\mathcal{U}_{n}(g,\beta)$
and $k=r_{+}+j$, then
\[
S_{k}\psi(\omega)=H_{+}+S_{j}\psi(\nu)\ge H_{+}+\beta j-K_{\beta}>0.
\]
The finitely many prefixes inside $\tau^{+}$ are non-negative as
well. Moreover
\[
S_{r_{+}+j}\psi(\omega)-\beta(r_{+}+j)=(H_{+}-\beta r_{+})+(S_{j}\psi(\nu)-\beta j),
\]
so one may choose a fixed $K_{{\rm U}}$. The fixed prefix $\tau^{+}$
changes the prefix estimate from \eqref{eq:bounded-window-frostman}
only by a multiplicative constant; the finitely many prefix lengths
$k\le r_{+}$ are absorbed into the same constant. This gives \eqref{eq:excursion-prefix-law}.

For recovery, write $S_{k}\psi=-\gamma k+e_{k}$ with $|e_{k}|\le K_{{\rm D}}\coloneqq K_{-\gamma}$.
If $H\ge\gamma n+K_{{\rm D}}$, then
\[
H+S_{k}\psi\ge\gamma(n-k)\ge0\qquad(0\le k\le n),
\]
which proves the recovery assertion. Finally, taking $k=n$ in the
three prefix estimates and summing the corresponding finite probability
laws gives the three lower partition-sum estimates. This completes
the proof.
\end{proof}

Only interior values $v\in(a_{-},a_{+})$ are used below. The endpoints
$a_{\pm}$, when needed, are obtained by approximation from the interior.

\section{The pressure $\cP^{\N}_{\alpha}$ and recurrent $\alpha$-escaping
sets}

\label{sec:general-reflected-pressure}

Throughout this section the alphabet is finite, $\psi$ is constant
on one-cylinders (in particular, the step set is finite), $\phi$
is Hölder, and we assume Lindley reflection. Other finite boundary
modifications would require the corresponding bounded communication
and bridge estimates and are not treated here.

We first express the $\alpha$-escaping pressure of the $\N$-extension
directly in terms of the pressure of the base system. We then use
this representation to choose the positive excursion and negative
recovery rates needed for the recurrent Moran construction.

\subsection{The pressure $\cP^{\N}_{\alpha}(g)$ from the base pressure}

For a Hölder potential~$g$, by~\Cref{prop:zero-speed-pressure},
\[
\cP^{\N}_{0}(g)=\inf_{r\le0}\cP^{0}(g+r\psi).
\]
Define the truncated base-pressure function 
\[
G_{g}(q)\coloneqq\max\left\{ \cP^{0}(g+q\psi),\cP^{\N}_{0}(g)\right\} =\max\left\{ \cP^{0}(g+q\psi),\inf_{r\le0}\cP^{0}(g+r\psi)\right\} .
\]
Since the step cocycle takes both positive and negative values, $q\mapsto\cP^{0}(g+q\psi)$
tends to $+\infty$ as $|q|\to\infty$ and is strictly convex. Its
concave Legendre transform is the $\alpha$-escape pressure of the
$\Z$-extension: 
\begin{equation}
\inf_{q\in\R}\left\{ \cP^{0}(g+q\psi)-q\alpha\right\} =\inf_{q\in\R}\cP^{0}\bigl(g+q(\psi-\alpha)\bigr)=\cP^{\Z}_{\alpha}(g).\label{eq:F-legendre-translation}
\end{equation}
In our previous paper \cite{GJK-realline}, this pressure was denoted
by $\cP(g,\psi-\alpha)$. We also write 
\[
\widehat{G}_{g}(\alpha)\coloneqq\inf_{q\in\R}\left\{ G_{g}(q)-q\alpha\right\} .
\]

\begin{lem}[Excursion decomposition for the $\N$-extension]
\label{lem:reflected-pressure-decomposition} For every Hölder potential
$g$ and every $\alpha\in\mathcal{A}_{+}\setminus\left\{ 0\right\} $,
\begin{equation}
\cP^{\N}_{\alpha}(g)=\sup_{\beta\in[\alpha,a_{+}]}\left\{ \left(1-\frac{\alpha}{\beta}\right)\cP^{\N}_{0}(g)+\frac{\alpha}{\beta}\cP^{\Z}_{\beta}(g)\right\} .\label{eq:marked-pressure-without-A}
\end{equation}
\end{lem}

\begin{proof}
Fix an endpoint tolerance $K$. Let a word of the $\N$-extension
of length $n$ end at height $\alpha n+O(K)$, and let $\ell$ be
its last visit to the fixed boundary region on which reflection can
act. Put $u=n-\ell$. Since $\alpha>0$ and the positive step size
is at most $a_{+}$, one has $u\ge(\alpha/a_{+})n+O(1)$. After changing
the cut by a uniformly bounded number of symbols, reflection is inactive
on the terminal block. If $d$ denotes its displacement in the $\Z$-extension
and $\beta=d/u$, then
\[
\beta\in[\alpha-o(1),a_{+}],\qquad\frac{u}{n}=\frac{\alpha}{\beta}+o(1),\qquad\frac{\ell}{n}=1-\frac{\alpha}{\beta}+o(1).
\]
The initial block ends in a fixed bounded window around height $0$
and therefore has exponential rate at most $\cP^{\N}_{0}(g)$. The
terminal block is a word of the $\Z$-extension with the indicated
endpoint displacement and has rate at most $\cP^{\Z}_{\beta}(g)$.
There are only $O(n^{2})$ choices of the cut time and the integer
endpoint displacement. A finite mesh in $\beta$, bounded distortion
and continuity of $\beta\mapsto\cP^{\Z}_{\beta}(g)$ therefore give
\[
\cP^{\N}_{\alpha}(g)\le\sup_{\beta\in[\alpha,a_{+}]}\left\{ \left(1-\frac{\alpha}{\beta}\right)\cP^{\N}_{0}(g)+\frac{\alpha}{\beta}\cP^{\Z}_{\beta}(g)\right\} .
\]

For the reverse inequality fix $\beta\in(\alpha,a_{+})$ and $\eta>0$.
By the definition of $\cP^{\N}_{0}(g)$, first choose a fixed endpoint
window and then arbitrarily large lengths $\ell$ for which the corresponding
partition sum $Z^{\N}_{\ell,K}(g,0)$ is at least $\exp\{\ell(\cP^{\N}_{0}(g)-\eta)\}$.
For every such $\ell$ choose an integer $u=u(\ell)$ with
\begin{equation}
u=\frac{\alpha}{\beta-\alpha}\,\ell+O(1).\label{eq:decomp-excursion-length}
\end{equation}
Then $\beta u-\alpha(\ell+u)=O(1)$. Fix any $\gamma_{0}\in(0,-a_{-})$.
By \Cref{prop:bounded-window-phase-families} applied to $(g,\eta,\beta,\gamma_{0})$,
for every sufficiently large $u$,
\[
\sum_{\omega\in\mathcal{U}_{u}(g,\beta)}\exp(\sup_{[\omega]}S_{u}g)\ge C^{-1}\exp\{u(\cP^{\Z}_{\beta}(g)-\eta)\},
\]
and every word in $\mathcal{U}_{u}(g,\beta)$ is a positive direct
route with $|S_{k}\psi-\beta k|\le K_{{\rm U}}$ for all $k\le u$.
Concatenate each word counted in $Z^{\N}_{\ell,K}(g,0)$ with each
such excursion word. The first block ends at bounded height, the excursion
word has non-negative partial displacement, and the concatenated endpoint
therefore satisfies
\[
h_{\ell+u}=\alpha(\ell+u)+O(1).
\]
Because the base is the full shift, no symbolic connector is required;
bounded distortion changes the product weight only by a uniform factor.
Dividing by $n=\ell+u$ and using \eqref{eq:decomp-excursion-length}
gives
\[
\cP^{\N}_{\alpha}(g)\ge\left(1-\frac{\alpha}{\beta}\right)\cP^{\N}_{0}(g)+\frac{\alpha}{\beta}\cP^{\Z}_{\beta}(g)-\eta.
\]
Letting $\eta\downarrow0$ and then taking the supremum over interior
$\beta$ proves the lower bound for $0<\alpha<a_{+}$; the values
$\beta=\alpha$ and $\beta=a_{+}$ follow by continuity. If $\alpha=a_{+}$,
a bounded endpoint defect permits only $O(1)$ non-maximal increments,
so the endpoint pressures of the $\N$- and $\Z$-extensions agree.
This completes the proof.
\end{proof}

To state the formula uniformly, set
\[
F_{g}(q)\coloneqq\cP^{0}(g+q\psi),
\]
and let $q_{\min}(g)$ be its unique minimizer. If $q_{\min}(g)=0$,
set $q^{*}(g)\coloneqq0$. Otherwise, let $q^{*}(g)$ be the unique
point on the opposite side of $q_{\min}(g)$ from $0$ such that
\begin{equation}
\cP^{0}(g+q^{*}(g)\psi)=\cP^{0}(g).\label{eq:pressure-partner}
\end{equation}
Strict convexity and $F_{g}(q)\to+\infty$ as $|q|\to\infty$ give
existence and uniqueness. Put
\[
a^{*}(g)\coloneqq\left.\partial_{q}\cP^{0}(g+q\psi)\right|_{q=q^{*}(g)}.
\]
Thus $a^{*}(g)>0$ exactly when $\left.\partial_{q}\cP^{0}(g+q\psi)\right|_{q=0}<0$.
\begin{thm}[Escape pressure of the $\N$-extension]
\label{thm:truncated-base-pressure} For every Hölder potential $g$
and every $\alpha\in\mathcal{A}_{+}$,
\begin{equation}
{\cP^{\N}_{\alpha}(g)=\widehat{G}_{g}(\alpha)=\inf_{q\in\R}\left\{ G_{g}(q)-q\alpha\right\} .}\label{eq:truncated-base-pressure}
\end{equation}
Moreover,
\begin{equation}
{\cP^{\N}_{\alpha}(g)=\begin{cases}
\cP^{0}(g)-q^{*}(g)\alpha, & 0\le\alpha<a^{*}(g),\\[1mm]
\cP^{\Z}_{\alpha}(g), & \max\{0,a^{*}(g)\}\le\alpha\le a_{+}.
\end{cases}}\label{eq:contact-formula}
\end{equation}
The first branch is void when $a^{*}(g)\le0$; when $a^{*}(g)>0$,
the two branches meet tangentially at $\alpha=a^{*}(g)$.
\end{thm}

\begin{proof}
At $\alpha=0$,
\[
\widehat{G}_{g}(0)=\inf_{q\in\R}G_{g}(q)=\cP^{\N}_{0}(g).
\]
For $\alpha>0$ write $c\coloneqq\cP^{\N}_{0}(g)$ and
\[
E_{\alpha}(\beta)\coloneqq\Bigl(1-\frac{\alpha}{\beta}\Bigr)c+\frac{\alpha}{\beta}\cP^{\Z}_{\beta}(g),\qquad\beta\in[\alpha,a_{+}],
\]
for the expression in braces in~\eqref{eq:marked-pressure-without-A}.
By~\Cref{lem:reflected-pressure-decomposition} we have to show that
$\sup_{\beta\in[\alpha,a_{+}]}E_{\alpha}(\beta)=\widehat{G}_{g}(\alpha)$.
For the inequality $\le$ fix $q\in\R$ and $\beta\in[\alpha,a_{+}]$.
By~\eqref{eq:F-legendre-translation}, $\cP^{\Z}_{\beta}(g)\le F_{g}(q)-q\beta$,
and hence
\[
E_{\alpha}(\beta)\le\Bigl(1-\frac{\alpha}{\beta}\Bigr)c+\frac{\alpha}{\beta}F_{g}(q)-q\alpha\le\max\{c,F_{g}(q)\}-q\alpha=G_{g}(q)-q\alpha.
\]
Taking the infimum over $q$ gives $E_{\alpha}(\beta)\le\widehat{G}_{g}(\alpha)$.
The reverse inequality follows from the explicit description of $\widehat{G}_{g}$
obtained in the remaining part of the proof: if $q^{*}(g)\le0$, or
if $q^{*}(g)>0$ and $\alpha\ge a^{*}(g)$, then $\widehat{G}_{g}(\alpha)=\cP^{\Z}_{\alpha}(g)=E_{\alpha}(\alpha)$;
if $q^{*}(g)>0$ and $\alpha<a^{*}(g)$, then $c=\cP^{0}(g)$, $a^{*}(g)\in(\alpha,a_{+})$
and, since the infimum defining $\cP^{\Z}_{a^{*}(g)}(g)$ is attained
at $q=q^{*}(g)$,
\[
\begin{aligned}E_{\alpha}(a^{*}(g)) & =\Bigl(1-\frac{\alpha}{a^{*}(g)}\Bigr)\cP^{0}(g)+\frac{\alpha}{a^{*}(g)}\bigl(\cP^{0}(g)-q^{*}(g)a^{*}(g)\bigr)\\
 & =\cP^{0}(g)-q^{*}(g)\alpha=\widehat{G}_{g}(\alpha).
\end{aligned}
\]
This proves~\eqref{eq:truncated-base-pressure}.

It remains to read off this transform from the geometry of $F_{g}$.
If $q^{*}(g)\le0$, then $q_{\min}(g)\le0$ and hence
\[
\cP^{\N}_{0}(g)=\min_{q\in\R}F_{g}(q),\qquad G_{g}=F_{g}.
\]
Thus $\cP^{\N}_{\alpha}(g)=\cP^{\Z}_{\alpha}(g)$ for every $\alpha\ge0$.
Since then $a^{*}(g)\le0$, this is precisely \eqref{eq:contact-formula}.

If $q^{*}(g)>0$, then $q_{\min}(g)>0$ and therefore $\cP^{\N}_{0}(g)=\cP^{0}(g)$.
By \eqref{eq:pressure-partner} and strict convexity,
\[
G_{g}(q)=\cP^{0}(g)\quad(0\le q\le q^{*}(g)),
\]
while $G_{g}(q)=F_{g}(q)$ outside this interval. Hence the minimizer
in $\inf_{q}\{G_{g}(q)-q\alpha\}$ is $q^{*}(g)$ for $0\le\alpha<a^{*}(g)$,
and lies in $[q^{*}(g),\infty)$ for $a^{*}(g)\le\alpha\le a_{+}$.
This gives \eqref{eq:contact-formula}. Finally, Legendre duality
gives
\[
\left.\partial_{\alpha}\cP^{\Z}_{\alpha}(g)\right|_{\alpha=a^{*}(g)}=-q^{*}(g),
\]
which is also the slope of the linear branch, proving tangency.
\end{proof}

\begin{cor}[Comparison of $\cP^{\N}_{\alpha}(g)$ and $\cP^{\Z}_{\alpha}(g)$]
\label{cor:reflected-pressure-dominates-translation} For every Hölder
potential $g$ and every $\alpha\in\mathcal{A}_{+}$,
\begin{align*}
\cP^{\N}_{\alpha}(g) & =\inf_{q\in\R}\left\{ \max\left(\cP^{0}\bigl(g+q(\psi-\alpha)\bigr),\inf_{r\le0}\cP^{0}(g+r\psi-q\alpha)\right)\right\} \\
 & \ge\inf_{q\in\R}\cP^{0}\bigl(g+q(\psi-\alpha)\bigr)=\cP^{\Z}_{\alpha}(g).
\end{align*}
\end{cor}

\begin{proof}
The identity is \Cref{thm:truncated-base-pressure} rewritten using
$G_{g}$ and the invariance of pressure under addition of constants.
The inequality follows from $\max\{A,B\}\ge A$, and the final identity
is \eqref{eq:F-legendre-translation}.
\end{proof}

\subsection{Geometric escape pressure}

For the geometric potential we write 
\[
p(s,q)\coloneqq\cP^{0}(s\phi+q\psi),\quad c_{s}\coloneqq p^{\N}_{0}(s)=\inf_{q\le0}p(s,q),\quad G_{s}(q)\coloneqq\max\{c_{s},p(s,q)\}.
\]
Thus \eqref{eq:truncated-base-pressure} becomes 
\[
p^{\N}_{\alpha}(s)=\inf_{q\in\R}\left\{ G_{s}(q)-q\alpha\right\} =\inf_{q\in\R}\left\{ \max\{p^{\N}_{0}(s),p(s,q)\}-q\alpha\right\} .
\]
Taking $g=s\phi$ in \Cref{cor:reflected-pressure-dominates-translation}
gives 
\[
p^{\N}_{\alpha}(s)\ge p^{\Z}_{\alpha}(s)\qquad(\alpha\in\mathcal{A}_{+}).
\]
Let $\mu_{s\phi}$ denote the equilibrium state of $s\phi$, that
is, the Gibbs measure with parameters $(s,q)=(s,0)$. Its drift, which
we also call the \emph{weighted drift} at $s$, is
\[
d(s)\coloneqq\partial_{q}p(s,0)=\int\psi\,\dd\mu_{s\phi}.
\]
 Recall that $\cP^{0}(\delta_{0}\phi)=0$ and $\partial_{q}p(\delta_{0},0)=\alpha_{\max}$.
Since the base is a finite full shift and $\phi$ is Hölder, the map
$(s,q)\longmapsto p(s,q)$ is real analytic. Moreover, the standing
assumption $a_{-}<0<a_{+}$ implies that $\psi$ is not cohomologous
to a constant, and hence $\partial_{qq}p(s,q)>0.$

For the geometric potential, write
\[
q^{*}(s)\coloneqq q^{*}(s\phi),\qquad a^{*}(s)\coloneqq\partial_{q}p(s,q^{*}(s)).
\]
By construction,
\[
a^{*}(s)>0\quad\Longleftrightarrow\quad\partial_{q}p(s,0)<0,
\]
so the single parameter $a^{*}(s)$ encodes the two regimes in \Cref{thm:truncated-base-pressure}.

\subsection{Poincaré exponent consequences}

By \Cref{prop:delta-pressure-root} and the preceding pressure inequality,
\[
\delta^{\N}(\alpha)\ge\delta^{\N}_{u}(\alpha).
\]

The contact formula also gives an intrinsic criterion for when the
endpoint exponent is strictly larger than the uniform exponent. It
involves only $\delta^{\N}_{u}(\alpha)$ and quantities associated
with the base pressure at this parameter.
\begin{cor}[Criterion for strict gain]
\label{cor:strict-reflection-gain} Let $\alpha\in(0,a_{+})$. Then
\[
\delta^{\N}(\alpha)>\delta^{\N}_{u}(\alpha)\quad\Longleftrightarrow\quad\alpha<a^{*}(\delta^{\N}_{u}(\alpha)).
\]
Moreover, for~$0<\alpha<a_{+}$,
\begin{equation}
\alpha\in\mathcal{A}^{\N}\quad\Longrightarrow\quad\delta^{\N}(\alpha)<\delta^{\N}(0),\label{eq:strict-gain-below-zero-threshold}
\end{equation}
or equivalently,
\begin{equation}
\delta^{\N}(\alpha)\ge\delta^{\N}(0)\quad\Longrightarrow\quad\delta^{\N}(\alpha)=\delta^{\N}_{u}(\alpha)\qquad(0<\alpha<a_{+}).\label{eq:large-regime-equality}
\end{equation}
\end{cor}

\begin{proof}
By \Cref{prop:uniform-exponent-comparison}, for $\alpha>0$ we have
$\delta^{\N}_{u}(\alpha)=\delta^{\Z}(\alpha)$. Put $s=\delta^{\N}_{u}(\alpha)$.
Then $p^{\Z}_{\alpha}(s)=0$, and strict monotonicity in the geometric
parameter gives
\[
\delta^{\N}(\alpha)>\delta^{\N}_{u}(\alpha)\quad\Longleftrightarrow\quad p^{\N}_{\alpha}(s)>0.
\]
By \Cref{thm:truncated-base-pressure}, if $a^{*}(s)\le0$ then $p^{\N}_{\alpha}(s)=p^{\Z}_{\alpha}(s)=0$
for every $\alpha\ge0$. If $a^{*}(s)>0$, the linear branch is tangent
to the strictly concave function $\eta\mapsto p^{\Z}_{\eta}(s)$ at
$\eta=a^{*}(s)$ and lies strictly above it for $0\le\eta<a^{*}(s)$.
Since $p^{\Z}_{\alpha}(s)=0$, this yields
\[
p^{\N}_{\alpha}(s)>0\quad\Longleftrightarrow\quad\alpha<a^{*}(s),
\]
which proves the criterion. For the additional statement, suppose~$\alpha\in\mathcal{A}^{\N}\cap(0,a_{+})$
and put~$s_{\N}\coloneqq\delta^{\N}(\alpha)$. Then $p^{\N}_{\alpha}(s_{\N})=0>p^{\Z}_{\alpha}(s_{\N})$,
so the contact branch in~\Cref{thm:truncated-base-pressure} is active
at~$s_{\N}$. Hence
\[
\alpha<a^{*}(s_{\N}),\qquad q^{*}(s_{\N})>0,\qquad0=p(s_{\N},0)-\alpha q^{*}(s_{\N}).
\]
Thus~$p^{\N}_{0}(s_{\N})=p(s_{\N},0)=\alpha q^{*}(s_{\N})>0$, and~\Cref{prop:delta-pressure-root}
gives~$s_{\N}<\delta^{\N}(0)$. This proves~\eqref{eq:strict-gain-below-zero-threshold};
the equivalent implication~\eqref{eq:large-regime-equality} follows
from~$\delta^{\N}\ge\delta^{\N}_{u}$.
\end{proof}

Thus strict gain can be decided without first computing $\delta^{\N}(\alpha)$:
one computes the intrinsic uniform exponent $\delta^{\N}_{u}(\alpha)$
and then checks the single contact escape rate $a^{*}(\delta^{\N}_{u}(\alpha))$.
Recall that
\[
\mathcal{A}^{\N}\coloneqq\{\alpha\in(0,a_{+}):\delta^{\N}(\alpha)>\delta^{\N}_{u}(\alpha)\}.
\]
Then \Cref{cor:strict-reflection-gain} gives the explicit description
displayed in the introduction.
\begin{rem}[No strict gain near the maximal escape rate]
\label{rem:no-gain-near-max-speed}There exists $\varepsilon>0$
such that 
\[
\delta^{\N}(\alpha)=\delta^{\N}_{u}(\alpha)\qquad(a_{+}-\varepsilon<\alpha\le a_{+}).
\]
 Indeed, the deflated pressure $R(s,q):=\int^{1}_{0}\partial_{q}p(s,tq)\,dt$,
which equals $(p(s,q)-p(s,0))/q$ for $q\ne0$ and is used again in
the proof of~\Cref{thm:regularity-reflected-spectrum}, has strictly
positive $q$-derivative, and its unique zero is $q^{*}(s)$. The
implicit-function theorem therefore makes $q^{*}$ and $a^{*}$ continuous,
including at parameters with $d(s)=0$. Set $s_{+}:=\delta^{\Z}(a_{+})$.
The equilibrium state of $s_{+}\phi+q^{*}(s_{+})\psi$ has full support,
so $a^{*}(s_{+})<a_{+}$. Continuity of $\delta^{\Z}$ at the extremal
escape rate $a_{+}$ (the map $\alpha\mapsto p^{\Z}_{\alpha}(s)$
is concave and upper semicontinuous on $[a_{-},a_{+}]$, hence continuous,
so its zero $\delta^{\Z}(\alpha)$ depends continuously on $\alpha\in[a_{-},a_{+}]$)
gives $a^{*}(\delta^{\Z}(\alpha))<\alpha$ for all $\alpha$ sufficiently
close to $a_{+}$. Apply \cref{cor:strict-reflection-gain} and the
endpoint identity $\delta^{\N}(a_{+})=\delta^{\N}_{u}(a_{+})$ of~\Cref{rem:endpoint-speeds}.
In particular, a strict-gain interval cannot extend all the way to
the maximal escape rate.
\end{rem}

\subsection{Regularity of the Poincaré spectrum of the $\N$-extension}

\begin{thm}[Regularity and geometry of the strict-gain set] \label{thm:regularity-reflected-spectrum}
Recall the strict-gain set 
\[
\mathcal{A}^{\N}\coloneqq\{\alpha\in(0,a_{+}):\delta^{\N}(\alpha)>\delta^{\N}_{u}(\alpha)\}.
\]
Then $\mathcal{A}^{\N}$ is relatively open in $(0,a_{+})$ and is
a finite union of intervals. Moreover, $\mathcal{A}^{\N}\subset(0,a_{+}-\varepsilon)$
for some $\varepsilon>0$ by~\Cref{rem:no-gain-near-max-speed}.
Consequently, its complement in $(0,a_{+})$ is a finite union of
relatively closed intervals, some of which may degenerate to isolated
points.

On every connected component of $\mathcal{A}^{\N}\cap(0,a_{+})$,
the map $\alpha\mapsto\delta^{\N}(\alpha)$ is real analytic. On the
interior of $(0,a_{+})\setminus\mathcal{A}^{\N}$, one has $\delta^{\N}=\delta^{\N}_{u}$,
and hence $\delta^{\N}$ is real analytic there as well.

Let $\alpha_{0}\in\partial\mathcal{A}^{\N}\cap(0,a_{+})$ and set
\[
s_{0}\coloneqq\delta^{\N}(\alpha_{0})=\delta^{\N}_{u}(\alpha_{0}).
\]
Then $\partial_{q}p(s_{0},0)<0$, so $q^{*}(s)>0$ in a neighborhood
of $s_{0}$. Put 
\[
H(s)\coloneqq p(s,0)-q^{*}(s)a^{*}(s).
\]
Then $\alpha_{0}=a^{*}(s_{0})$, $H(s_{0})=0$ and $\delta^{\N}$
is $C^{1}$ at $\alpha_{0}$. If $\alpha_{0}$ is not an isolated
point of $(0,a_{+})\setminus\mathcal{A}^{\N}$, then 
\begin{equation}
\delta^{\N}\text{ is }C^{2}\text{ at }\alpha_{0}\quad\Longleftrightarrow\quad H'(s_{0})=0.\label{eq:C2-delta-contact-criterion}
\end{equation}
Equivalently, $1-a^{*}{}'(s_{0})(\delta^{\N}_{u})'(\alpha_{0})=0,$
or $\partial_{s}p(s_{0},q^{*}(s_{0}))-q^{*}(s_{0})a^{*}{}'(s_{0})=0.$
In particular, a simple zero of $H$ gives a genuine second-order
phase transition: $\delta^{\N}$ is $C^{1}$ but not $C^{2}$. If
instead $\alpha_{0}$ is an isolated point of $(0,a_{+})\setminus\mathcal{A}^{\N}$,
then $\delta^{\N}$ is real analytic in a neighborhood of $\alpha_{0}$.

Finally, 
\[
\alpha_{\max}<0\quad\Longleftrightarrow\quad\exists\,\varepsilon>0:\ (0,\varepsilon)\subset\mathcal{A}^{\N},
\]
which characterizes the jump of the uniform spectrum at $\alpha=0$;
see~\Cref{cor:continuity-zero-jump}. \end{thm}
\begin{proof}
\emph{Continuity and analyticity away from transition points.} By
\Cref{prop:delta-pressure-root}, $\delta^{\N}(\alpha)$ is the unique
zero of $s\mapsto p^{\N}_{\alpha}(s)$. Since $\phi\le-\kappa<0$,
for $t>0$ one has 
\[
p^{\N}_{\alpha}(s+t)\le p^{\N}_{\alpha}(s)-\kappa t.
\]
Thus the zero is simple, uniformly in $\alpha$. The same facts hold
for $p^{\Z}_{\alpha}$ and $\delta^{\N}_{u}$.

For fixed $s$, the Legendre formulas show that $\alpha\mapsto p^{\N}_{\alpha}(s)$
and $\alpha\mapsto p^{\Z}_{\alpha}(s)$ are finite and concave. They
are therefore continuous on $(0,a_{+})$; at $0$, concavity together
with upper semicontinuity gives right-continuity. The endpoint exponent
$\delta^{\N}$ is continuous on $[0,a_{+})$, while $\delta^{\N}_{u}$
is real analytic on $(0,a_{+})$ by \Cref{prop:uniform-exponent-comparison}
and \Cref{thm:Z-spectrum}. Hence $\mathcal{A}^{\N}$ is relatively
open in $(0,a_{+})$.

Let $\alpha\in\mathcal{A}^{\N}\cap(0,a_{+})$ and put $s=\delta^{\N}(\alpha)$.
Since $s>\delta^{\N}_{u}(\alpha)$, strict decrease gives $p^{\Z}_{\alpha}(s)<0=p^{\N}_{\alpha}(s)$.
Thus the contact branch is active and its inequality is strict: 
\[
p(s,0)-\alpha q^{*}(s)=0,\qquad\alpha<a^{*}(s).
\]
In particular, $\partial_{q}p(s,0)<0$. The contact equation $p(s,q^{*}(s))=p(s,0)$,
together with $\partial_{q}p(s,q^{*}(s))=a^{*}(s)>0$, shows by the
implicit-function theorem that $q^{*}$ is real analytic near $s$.
In the strict-contact region, $p^{\N}_{\alpha}(s)=p(s,0)-\alpha q^{*}(s)$;
its $s$-derivative does not vanish because this pressure is strictly
decreasing. A second application of the implicit-function theorem
therefore proves that $\delta^{\N}$ is real analytic on every component
of $\mathcal{A}^{\N}\cap(0,a_{+})$. On the interior of the complement,
$\delta^{\N}=\delta^{\N}_{u}$, so analyticity follows from \Cref{prop:uniform-exponent-comparison}
and \Cref{thm:Z-spectrum}.

\emph{Finiteness of the transition set.} Recall that $d(s)=\partial_{q}p(s,0)$.
For $q\ne0$, define 
\[
R(s,q)\coloneqq\frac{p(s,q)-p(s,0)}{q},
\]
and extend it analytically to $q=0$ by $R(s,0)=d(s)$. Strict convexity
and $a_{-}<0<a_{+}$ imply that $R(s,\cdot)$ has a unique zero $\bar{q}(s)$.
It is positive when $d(s)<0$, zero when $d(s)=0$, and negative when
$d(s)>0$. At a nonzero zero, 
\[
\partial_{q}R(s,\bar{q}(s))=\frac{\partial_{q}p(s,\bar{q}(s))}{\bar{q}(s)}>0,
\]
whereas $\partial_{q}R(s,0)=\tfrac{1}{2}\partial_{qq}p(s,0)>0$. Hence
$\bar{q}$ is real analytic on a neighborhood of $[0,\delta_{0}]$.
Moreover, $\bar{q}$ is the analytic continuation of $q^{*}$, so
$a^{*}(s)=\partial_{q}p(s,\bar{q}(s))$ is real analytic on the same
neighborhood. Consequently, 
\[
\bar{H}(s)\coloneqq p(s,0)-\bar{q}(s)\partial_{q}p(s,\bar{q}(s))
\]
is real analytic there. It is not identically zero. Indeed, if $d(\delta_{0})\ne0$,
then 
\[
\bar{H}(\delta_{0})=-\bar{q}(\delta_{0})\partial_{q}p(\delta_{0},\bar{q}(\delta_{0}))<0;
\]
if $d(\delta_{0})=0$, then $\bar{H}'(\delta_{0})=\partial_{s}p(\delta_{0},0)<0$.
Thus $\bar{H}$ has only finitely many zeros in $[0,\delta_{0}]$.

Now let $\alpha_{0}\in\partial\mathcal{A}^{\N}\cap(0,a_{+})$ and
put $s_{0}=\delta^{\N}_{u}(\alpha_{0})$. By the strict-gain criterion,
\[
\mathcal{A}^{\N}\cap(0,a_{+})=\{\alpha\in(0,a_{+}):\alpha<a^{*}(\delta^{\N}_{u}(\alpha))\}.
\]
The function $\alpha\mapsto a^{*}(\delta^{\N}_{u}(\alpha))-\alpha$
is continuous. Therefore equality holds at a boundary point: 
\[
\alpha_{0}=a^{*}(s_{0}).
\]
Since $\alpha_{0}>0$, the relation $a^{*}(s)>0\Longleftrightarrow d(s)<0$
gives $d(s_{0})<0$, and hence $\bar{q}(s_{0})=q^{*}(s_{0})>0$. Moreover,
\[
0=p^{\Z}_{\alpha_{0}}(s_{0})=p^{\Z}_{a^{*}(s_{0})}(s_{0})=H(s_{0})=\bar{H}(s_{0}).
\]
Thus every positive boundary point determines a zero of $\bar{H}$;
conversely, each such zero determines at most one boundary point through
$\alpha=a^{*}(s)$. The positive boundary is therefore finite. Together
with relative openness, this proves that $\mathcal{A}^{\N}$ is a
finite union of intervals.

\emph{Matching the analytic branches.} Fix a positive boundary point
$\alpha_{0}$. Since $d(s_{0})<0$, the contact zero $\delta_{c}(\alpha)$
is defined analytically near $\alpha_{0}$ by 
\[
p(\delta_{c}(\alpha),0)-\alpha q^{*}(\delta_{c}(\alpha))=0.
\]
At $\alpha_{0}$, it satisfies $\delta_{c}(\alpha_{0})=s_{0}=\delta^{\N}_{u}(\alpha_{0})$.
Differentiating the contact identity $p(s,q^{*}(s))=p(s,0)$ gives
\[
\partial_{s}p(s,q^{*}(s))+a^{*}(s)q^{*}{}'(s)=\partial_{s}p(s,0).
\]
Hence, at $(s_{0},\alpha_{0})$, 
\[
\partial_{s}\bigl(p(s,0)-\alpha q^{*}(s)\bigr)=\partial_{s}p(s_{0},q^{*}(s_{0}))<0.
\]
On the uniform branch, the minimizing parameter at contact is the
same $q^{*}(s_{0})$. Implicit differentiation of the two zero equations
yields 
\[
\delta_{c}'(\alpha_{0})=(\delta^{\N}_{u})'(\alpha_{0})=\frac{q^{*}(s_{0})}{\partial_{s}p(s_{0},q^{*}(s_{0}))}.
\]
Thus the two analytic branches have first-order contact, and $\delta^{\N}$
is $C^{1}$ at $\alpha_{0}$.

To compare second derivatives, write 
\[
C(s,\alpha)\coloneqq p(s,0)-\alpha q^{*}(s),\qquad T(s,\alpha)\coloneqq p^{\Z}_{\alpha}(s).
\]
For every nearby $s$, the functions $C$ and $T$, as well as their
first $\alpha$-derivatives, agree at $\alpha=a^{*}(s)$, whereas
$\partial_{\alpha\alpha}C(s,a^{*}(s))=0$, $\partial_{\alpha\alpha}T(s,a^{*}(s))=-1/\partial_{qq}p(s,q^{*}(s)).$
Analytic division therefore gives 
\begin{equation}
C(s,\alpha)-T(s,\alpha)=(\alpha-a^{*}(s))^{2}R_{*}(s,\alpha),\qquad R_{*}(s_{0},\alpha_{0})=\frac{1}{2\partial_{qq}p(s_{0},q^{*}(s_{0}))}>0.\label{eq:pressure-contact-square-factor}
\end{equation}
Evaluating this identity on the zero curve of $T$ and using the implicit-function
theorem for the zero of $C$ gives 
\begin{equation}
\delta_{c}''(\alpha_{0})-(\delta^{\N}_{u})''(\alpha_{0})=-\frac{\bigl(1-a^{*}{}'(s_{0})(\delta^{\N}_{u})'(\alpha_{0})\bigr)^{2}}{\partial_{s}p(s_{0},q^{*}(s_{0}))\,\partial_{qq}p(s_{0},q^{*}(s_{0}))}.\label{eq:delta-second-derivative-jump}
\end{equation}
Finally, $H'(s_{0})=\partial_{s}p(s_{0},q^{*}(s_{0}))-q^{*}(s_{0})a^{*}{}'(s_{0}),$
and the first-derivative formula gives 
\[
1-a^{*}{}'(s_{0})(\delta^{\N}_{u})'(\alpha_{0})=\frac{H'(s_{0})}{\partial_{s}p(s_{0},q^{*}(s_{0}))}.
\]
If $\alpha_{0}$ is not isolated in the complement of $\mathcal{A}^{\N}$,
the i.o.\ spectrum switches between the two analytic branches at
$\alpha_{0}$. Equations \eqref{eq:delta-second-derivative-jump}
and the preceding identity show that their second derivatives agree
exactly when $H'(s_{0})=0$, proving \eqref{eq:C2-delta-contact-criterion}.
If $\alpha_{0}$ is isolated in the complement, then the strict-gain
condition holds on both sides. Hence $\delta^{\N}=\delta_{c}$ on
a punctured neighborhood and, by continuity, at its center; thus $\delta^{\N}$
is real analytic there.

We are left to prove the equivalence concerning $\alpha=0$. Put $s_{0}\coloneqq\delta^{\Z}(0)$,
so that $s_{0}=\lim_{\alpha\downarrow0}\delta^{\N}_{u}(\alpha)$ by~\Cref{prop:uniform-exponent-comparison}
and continuity of $\delta^{\Z}$. Since $p^{\Z}_{0}(s)=\min_{q\in\R}p(s,q)\le p(s,0)=\cP^{0}(s\phi)$
and $s\mapsto p^{\Z}_{0}(s)$ is strictly decreasing with zero $s_{0}$,
we have $s_{0}\le\delta_{0}$, with equality if and only if $\min_{q}p(\delta_{0},q)$
is attained at $q=0$, that is, if and only if $\alpha_{\max}=d(\delta_{0})=0$.
Moreover, if $d(s_{1})=0$ for some $s_{1}\in[s_{0},\delta_{0}]$,
then $p^{\Z}_{0}(s_{1})=p(s_{1},0)=\cP^{0}(s_{1}\phi)\ge0=p^{\Z}_{0}(s_{0})$,
so strict monotonicity forces $s_{1}=s_{0}$, and then $\cP^{0}(s_{0}\phi)=0$
gives $s_{0}=\delta_{0}$. Hence $d$ has no zero in $[s_{0},\delta_{0}]$
unless $s_{0}=\delta_{0}$. Suppose first that $\alpha_{\max}<0$.
Then $s_{0}<\delta_{0}$ and $d(\delta_{0})<0$, so $d<0$ on $[s_{0},\delta_{0}]$
by the intermediate value theorem. In particular $a^{*}(s_{0})>0$,
and since $\alpha\mapsto a^{*}(\delta^{\N}_{u}(\alpha))-\alpha$ extends
continuously to $\alpha=0$ with value $a^{*}(s_{0})>0$, the strict-gain
criterion~\Cref{cor:strict-reflection-gain} gives $(0,\varepsilon)\subset\mathcal{A}^{\N}$
for some $\varepsilon>0$. Suppose next that $\alpha_{\max}>0$. Then
again $s_{0}<\delta_{0}$, and now $d>0$ on $[s_{0},\delta_{0}]$,
so $a^{*}(\delta^{\N}_{u}(\alpha))\le0<\alpha$ for all sufficiently
small $\alpha>0$, and $(0,\varepsilon)\cap\mathcal{A}^{\N}=\varnothing$
for some $\varepsilon>0$. Suppose finally that $\alpha_{\max}=0$.
Then $s_{0}=\delta_{0}$, $\bar{q}(\delta_{0})=0$ and $a^{*}(\delta_{0})=d(\delta_{0})=0$.
The function $\varrho(\alpha)\coloneqq a^{*}(\delta^{\Z}(\alpha))-\alpha$
is real analytic near $\alpha=0$ by~\Cref{thm:Z-spectrum} and the
analyticity of $a^{*}$, with $\varrho(0)=0$ and $\varrho'(0)=a^{*}{}'(\delta_{0})\,(\delta^{\Z})'(0)-1=-1$,
because $\delta^{\Z}$ attains its maximum at $\alpha_{\max}=0$.
Hence $\varrho<0$ on some interval $(0,\varepsilon)$, which is therefore
disjoint from $\mathcal{A}^{\N}$. This proves the equivalence.
\end{proof}

\begin{cor}[Contact multiplicity and transition order]
\label{cor:contact-multiplicity} Let $\alpha_{0}\in\partial\mathcal{A}^{\N}\cap(0,a_{+})$,
put $s_{0}:=\delta^{\N}_{u}(\alpha_{0})$, and let $m\ge1$ be the
multiplicity of the zero of $H$ at $s_{0}$ (recall that $H(s_{0})=0$
and that $H$ does not vanish identically, by the proof of~\Cref{thm:regularity-reflected-spectrum}).
For the analytic contact branch $\delta_{c}$ from the proof of \Cref{thm:regularity-reflected-spectrum},
there is a real-analytic function $W$, with $W(\alpha_{0})>0$, such
that 
\begin{equation}
\delta_{c}(\alpha)-\delta^{\N}_{u}(\alpha)=(\alpha-\alpha_{0})^{2m}W(\alpha)\label{eq:contact-multiplicity-factor}
\end{equation}
 near $\alpha_{0}$. If $\alpha_{0}$ is a genuine change of branches,
that is, if $\mathcal{A}^{\N}$ lies on exactly one side of $\alpha_{0}$,
then $m$ is odd and $\delta^{\N}$ is $C^{2m-1}$ but not $C^{2m}$
there. If $\alpha_{0}$ is isolated in the complement of $\mathcal{A}^{\N}$,
then $m$ is even and $\delta^{\N}$ is real analytic near $\alpha_{0}$.
\end{cor}

\begin{proof}
Use $T(s,\alpha)=p^{\Z}_{\alpha}(s)$ and $C(s,\alpha)=p(s,0)-\alpha q^{*}(s)$
as above. At contact, $T_{\alpha}=-q^{*}(s_{0})\ne0$, so the equation
$T(s,b(s))=0$ defines an analytic function $b$ near $s_{0}$, with
$b(s_{0})=\alpha_{0}$ and $b'(s_{0})=p_{s}(s_{0},q^{*}(s_{0}))/q^{*}(s_{0})\ne0$.
Since $T(s,a^{*}(s))=H(s)$, analytic division gives 
\[
H(s)=(a^{*}(s)-b(s))U(s),\qquad U(s_{0})=-q^{*}(s_{0})\ne0.
\]
 The functions $b$ and $\delta^{\N}_{u}$ are local inverses. Hence
$f(\alpha):=\alpha-a^{*}(\delta^{\N}_{u}(\alpha))$ has a zero of
precisely the same multiplicity $m$ at $\alpha_{0}$. Evaluating
\cref{eq:pressure-contact-square-factor} on the zero curve of $T$
and factoring $C$ along its analytic zero curve yields 
\[
\delta_{c}(\alpha)-\delta^{\N}_{u}(\alpha)=f(\alpha)^{2}V(\alpha),\qquad V(\alpha_{0})>0,
\]
 because $R_{*}(s_{0},\alpha_{0})>0$ and $C_{s}(s_{0},\alpha_{0})<0$.
Factoring the zero of $f$ proves \cref{eq:contact-multiplicity-factor}.

By \cref{cor:strict-reflection-gain}, strict gain is equivalent to
$f<0$. At a genuine change of branches, $f$ changes sign, so $m$
is odd. The two analytic branches then agree through derivative order
$2m-1$, but their derivatives of order $2m$ differ by $(2m)!W(\alpha_{0})>0$.
If the complement is isolated at $\alpha_{0}$, then $f<0$ on both
sides, so $m$ is even; the spectrum agrees with $\delta_{c}$ throughout
a neighborhood. These are the only possibilities because the boundary
is finite.
\end{proof}

\subsection{The recurrent $\alpha$-escaping set of the $\N$-extension}


We now turn to recurrent i.o.\ escape. The parameter choice is made
entirely from the pressure geometry above. Once the positive excursion
rate $\beta$ and recovery rate $\gamma$ have been chosen, the lower
bound is obtained from the bounded-window phase families of \Cref{prop:bounded-window-phase-families}.
\begin{thm}[Recurrent $\alpha$-escaping set of the $\N$-extension]
\label{thm:recurrent-reflected-escape} For every $\alpha\in(0,a_{+})$,
\begin{equation}
{\dim_{H}\bigl(E^{\N}(\alpha)\cap E^{\N}(0)\bigr)=\min\bigl\{\delta^{\N}(\alpha),\delta^{\N}(0)\bigr\}.}\label{eq:recurrent-reflected-dimension}
\end{equation}
\end{thm}

\begin{proof}
Since by \Cref{thm:uniform-N}, $\dim_{H}E^{\N}(0)=\delta^{\N}(0)$
and by a standard covering argument, $E^{\N}(\alpha)\leq\delta^{\N}(\alpha)$,
we get immediately get the upper bound
\[
\dim_{H}\bigl(E^{\N}(\alpha)\cap E^{\N}(0)\bigr)\le\min\left\{ \delta^{\N}(\alpha),\delta^{\N}(0)\right\} .
\]

For the reverse inequality fix $s<\min\{\delta^{\N}(\alpha),\delta^{\N}(0)\}.$
By \Cref{prop:delta-pressure-root} and strict monotonicity in the
geometric parameter,
\[
p^{\N}_{\alpha}(s)>0,\qquad c_{s}\coloneqq p^{\N}_{0}(s)>0.
\]

\medskip{}
\noindent\emph{Choice of phase rates.} By \eqref{eq:marked-pressure-without-A}
and $p^{\N}_{\alpha}(s)>0$, there is $\beta\in[\alpha,a_{+}]$ such
that
\[
\left(1-\frac{\alpha}{\beta}\right)c_{s}+\frac{\alpha}{\beta}p^{\Z}_{\beta}(s)>0.
\]
Because $0<\alpha<a_{+}$ and $\beta\mapsto p^{\Z}_{\beta}(s)$ is
continuous, we may perturb the choice, if necessary, and assume
\[
\alpha<\beta<a_{+},\qquad\left(1-\frac{\alpha}{\beta}\right)c_{s}+\frac{\alpha}{\beta}p^{\Z}_{\beta}(s)>0.
\]
For recovery use the minimizer $q^{\N}_{s}$ and the rate $v^{\N}_{s}$
from \eqref{eq:complementarity-waiting}. If $v^{\N}_{s}<0$, take
$\gamma\coloneqq-v^{\N}_{s}$. Then $p^{\Z}_{-\gamma}(s)=c_{s}>0.$
If $v^{\N}_{s}=0$, continuity of $v\mapsto p^{\Z}_{v}(s)$ and $p^{\Z}_{0}(s)=c_{s}>0$
give a sufficiently small $\gamma>0$ with $p^{\Z}_{-\gamma}(s)>0$.
In either case we can choose
\[
0<\gamma<-a_{-},\qquad p^{\Z}_{-\gamma}(s)>0.
\]
Choose $\eta>0$ so small that
\[
\begin{aligned}c_{s}-\eta & >0,\\
A_{\eta} & \coloneqq\left(1-\frac{\alpha}{\beta}\right)(c_{s}-\eta)+\frac{\alpha}{\beta}\bigl(p^{\Z}_{\beta}(s)-\eta\bigr)>0,\\
p^{\Z}_{-\gamma}(s)-\eta & >0.
\end{aligned}
\]
Fix the phase families from \Cref{prop:bounded-window-phase-families}
for these parameters.

\medskip{}
\noindent\emph{The Moran construction.} Set $M_{0}\coloneqq0$. Suppose
the $(j-1)$st recovery phase has ended at time $M_{j-1}$ in the
fixed boundary region supplied by \Cref{prop:bounded-window-phase-families}.
Choose $N_{j}>M_{j-1}$ so large that all phase lengths below are
in the range of that proposition and
\begin{equation}
\frac{M_{j-1}}{N_{j}}\longrightarrow0.\label{eq:moran-separation}
\end{equation}
No further scale separation is required.

Choose integers
\begin{align*}
u_{j} & =\frac{\alpha}{\beta}N_{j}+O(1),\\
\ell_{j} & =N_{j}-M_{j-1}-u_{j},\\
d_{j} & =\frac{\alpha}{\gamma}N_{j}+O(1).
\end{align*}
The bounded $O(1)$ in the choice of $d_{j}$ is chosen once and for
all so that every possible excursion endpoint, which is $\alpha N_{j}+O(1)$
uniformly over the family, lies at least $\gamma d_{j}+K_{{\rm D}}$
above the boundary, while its excess over $\gamma d_{j}$ remains
bounded. This is possible because all endpoint errors in \cref{prop:bounded-window-phase-families}
are uniform. From \eqref{eq:moran-separation},
\[
\frac{\ell_{j}}{N_{j}}\longrightarrow1-\frac{\alpha}{\beta},\qquad\frac{u_{j}}{N_{j}}\longrightarrow\frac{\alpha}{\beta}.
\]
In particular $\ell_{j},u_{j},d_{j}\to\infty$.

At generation $j$ choose independently a waiting word from $\mathcal{W}_{\ell_{j}}(s)$,
an excursion word from $\mathcal{U}_{u_{j}}(s,\beta)$ and a recovery
word from $\mathcal{D}_{d_{j}}(s,\gamma)$, with the corresponding
probability laws of \cref{prop:bounded-window-phase-families}. The
base is a full shift, so these words concatenate without symbolic
connectors.

The waiting word keeps the orbit in a bounded boundary strip. The
excursion is a positive direct route and
\[
S_{u_{j}}\psi=\beta u_{j}+O(1)=\alpha N_{j}+O(1).
\]
Since the excursion begins at a bounded height, at the marked time
$N_{j}$ we therefore have
\begin{equation}
h(F^{N_{j}}_{\N}x)=\alpha N_{j}+O(1),\label{eq:moran-marked-height}
\end{equation}
uniformly over all branches. By the choice of $d_{j}$, the recovery
word remains above the boundary for its whole duration and finishes
in a bounded region. Put
\[
M_{j}\coloneqq N_{j}+d_{j}.
\]
The product of all phase laws defines a probability measure $\nu$
on a Moran set $\mathcal{M}$. Equation \eqref{eq:moran-marked-height}
holds for infinitely many $j$ and the orbit returns to a bounded
region at every $M_{j}$. Hence
\begin{equation}
\mathcal{M}\subset E^{\N}(\alpha)\cap E^{\N}(0).\label{eq:moran-recurrent-inclusion}
\end{equation}

\medskip{}
\noindent\emph{The Frostman estimate.} Let $r$ be an arbitrary cylinder
length in generation $j$. Multiplying the prefix estimates in \cref{prop:bounded-window-phase-families}
and using the Bowen property of $s\phi$ across the at most three
phase joins in each completed generation gives
\begin{equation}
\log\nu([\omega|_{r}])\le\sup_{[\omega|_{r}]}S_{r}(s\phi)-\Pi_{\eta}(r)+O(j).\label{eq:moran-prefix}
\end{equation}
Here $\Pi_{\eta}(r)$ is obtained by assigning the rate $c_{s}-\eta$
to waiting time, $p^{\Z}_{\beta}(s)-\eta$ to excursion time and $p^{\Z}_{-\gamma}(s)-\eta$
to recovery time.

For all sufficiently large generations this accumulated increase in
pressure is non-negative at every point in time. During waiting it
increases with slope $c_{s}-\eta>0$. During the excursion the new-generation
contribution after $k$ excursion steps is
\[
\ell_{j}(c_{s}-\eta)+k\bigl(p^{\Z}_{\beta}(s)-\eta\bigr),\qquad0\le k\le u_{j}.
\]
This is affine in $k$. Its value at $k=0$ is non-negative, while
at $k=u_{j}$,
\[
\begin{aligned}\ell_{j}(c_{s}-\eta)+u_{j}\bigl(p^{\Z}_{\beta}(s)-\eta\bigr) & =N_{j}A_{\eta}-M_{j-1}(c_{s}-\eta)+O(1)\\
 & >0
\end{aligned}
\]
for large $j$ by \eqref{eq:moran-separation} and $A_{\eta}>0$.
Hence it is non-negative throughout the excursion. Recovery then increases
it further because $p^{\Z}_{-\gamma}(s)-\eta>0$. Inductively,
\begin{equation}
\Pi_{\eta}(r)\ge0\label{eq:moran-pressure-surplus}
\end{equation}
at every sufficiently large cylinder scale.

The scale separation \eqref{eq:moran-separation} implies $j=o(r)$
for cylinder lengths in generation $j$: eventually $N_{j}\ge2M_{j-1}$,
so the cycle endpoints grow at least geometrically. Combining \eqref{eq:moran-prefix}
and \eqref{eq:moran-pressure-surplus} therefore gives
\[
\log\nu([\omega|_{r}])\le\sup_{[\omega|_{r}]}S_{r}(s\phi)+o(r).
\]
Bounded distortion gives $\log|J_{\omega|_{r}}|=\sup_{[\omega|_{r}]}S_{r}\phi+O(1),$
and uniform expansion gives $\sup S_{r}\phi\le-\kappa r$ for some
$\kappa>0$. Hence, for every $s'<s$,
\[
\nu([\omega|_{r}])\le C_{s'}|J_{\omega|_{r}}|^{s'}
\]
for all sufficiently long cylinders meeting $\mathcal{M}$. A standard
stopping-cylinder argument for a finite full-branch expanding system
converts this cylinder estimate into the same estimate for intervals:
an interval of a given small size meets only a uniformly bounded number
of maximal cylinders whose diameter has just dropped below that size.
The mass distribution principle therefore yields
\[
\dim_{H}\mathcal{M}\ge s'.
\]
Letting $s'\uparrow s$, then $s\uparrow\min\{\delta^{\N}(\alpha),\delta^{\N}(0)\}$,
and using \eqref{eq:moran-recurrent-inclusion} proves the lower bound
and hence \eqref{eq:recurrent-reflected-dimension}.
\end{proof}


\subsection{The full $\alpha$-escaping set of the $\N$-extension}


The preceding theorem determines the recurrent contribution. The uniformly
$\alpha$-escaping and transient parts of the remaining case contribute
the same Hausdorff dimension.
\begin{cor}[Full $\alpha$-escaping set of the $\N$-extension]
\label{cor:full-reflected-escape} For every $\alpha\in(0,a_{+})$,
\[
\dim_{H}E^{\N}(\alpha)=\delta^{\N}(\alpha).
\]
\end{cor}

\begin{proof}
The upper bound follows from a standard covering argument. To prove
the lower bound, suppose first that $\delta^{\N}(\alpha)\le\delta^{\N}(0).$
Then \Cref{thm:recurrent-reflected-escape} gives 
\[
\dim_{H}\bigl(E^{\N}(\alpha)\cap E^{\N}(0)\bigr)=\delta^{\N}(\alpha).
\]
Since this set is contained in $E^{\N}(\alpha)$, the required lower
bound follows.

Now suppose that $\delta^{\N}(\alpha)>\delta^{\N}(0).$ By \Cref{cor:strict-reflection-gain},
$\delta^{\N}_{u}(\alpha)=\delta^{\N}(\alpha).$ Hence \Cref{thm:uniform-N}
gives 
\[
\dim_{H}E^{\N}_{u}(\alpha)=\dim_{H}\bigl(E^{\N}(\alpha)\setminus E^{\N}\left(0\right)\bigr)=\delta^{\N}(\alpha).
\]
Since $E^{\N}_{u}(\alpha)\subset E^{\N}(\alpha)$, this gives the
required lower bound for $\dim_{H}E^{\N}(\alpha)$.
\end{proof}

\begin{rem}[Endpoint escape rates]
\label{rem:endpoint-speeds}At $\alpha=0$, 
\[
\dim_{H}E^{\N}(0)=\delta^{\N}(0)
\]
by \Cref{thm:uniform-N}. At $\alpha=a_{+}$, the usual extremal-step
argument gives 
\[
p^{\N}_{a_{+}}(s)=p^{\Z}_{a_{+}}(s).
\]
Indeed, if 
\[
D_{n}\coloneqq a_{+}n-S^{\N}_{n}\psi,
\]
then $D_{n}$ is non-decreasing, and every increment strictly smaller
than $a_{+}$ increases $D_{n}$ by a fixed positive amount. Thus
\[
|S^{\N}_{n}\psi-a_{+}n|\le K
\]
allows only $O_{K}(1)$ non-extremal increments. Hence the endpoint
pressures of the $\N$- and $\Z$-extensions agree, and $\dim_{H}E^{\N}(a_{+})=\delta^{\N}(a_{+})=\delta^{\N}_{u}(a_{+}).$
Consequently, 
\[
\dim_{H}E^{\N}(\alpha)=\delta^{\N}(\alpha)\qquad(0\le\alpha\le a_{+}).
\]
\end{rem}

\begin{proof}[Proof of~\Cref{thm:exact-reflected-spectrum}]
 The upper bound in~\eqref{eq:exact-reflected-spectrum} follows
from a standard covering argument. The identity~\eqref{eq:exact-reflected-spectrum}
is~\Cref{cor:full-reflected-escape} for $0<\alpha<a_{+}$ and~\Cref{rem:endpoint-speeds}
for $\alpha\in\{0,a_{+}\}$, and~\eqref{eq:main-recurrent-dimension}
is~\Cref{thm:recurrent-reflected-escape}. For~\eqref{eq:recurrent-realizes}
suppose that $\delta^{\N}_{u}(\alpha)\le\delta^{\N}(0)$. If we had
$\delta^{\N}(\alpha)>\delta^{\N}(0)$, then~\eqref{eq:large-regime-equality}
would give $\delta^{\N}(\alpha)=\delta^{\N}_{u}(\alpha)\le\delta^{\N}(0)$,
a contradiction. Hence $\delta^{\N}(\alpha)\le\delta^{\N}(0)$, and~\eqref{eq:main-recurrent-dimension}
yields $\dim_{H}(E^{\N}(\alpha)\cap E^{\N}(0))=\delta^{\N}(\alpha)$.
For~\eqref{eq:translation-realizes} suppose that $\delta^{\N}_{u}(\alpha)>\delta^{\N}(0)$.
Since $\delta^{\N}\ge\delta^{\N}_{u}$ by~\Cref{cor:reflected-pressure-dominates-translation}
and~\Cref{prop:delta-pressure-root}, we get $\delta^{\N}(\alpha)\ge\delta^{\N}_{u}(\alpha)>\delta^{\N}(0)$,
so~\eqref{eq:large-regime-equality} gives $\delta^{\N}(\alpha)=\delta^{\N}_{u}(\alpha)$,
and~\Cref{thm:uniform-N} gives $\dim_{H}E^{\N}(\alpha)\setminus E^{\N}\left(0\right)=\delta^{\N}_{u}(\alpha)$.
\end{proof}

\begin{cor}[Continuity and the jump at $\alpha=0$]
\label{cor:continuity-zero-jump} The map~$\alpha\mapsto\dim_{H}E^{\N}(\alpha)$
is continuous on~$[0,a_{+})$. Moreover,~$\alpha\mapsto\dim_{H}E^{\N}_{u}(\alpha)$
is discontinuous at~$0$ if and only if
\[
\alpha_{\max}<0,
\]
and this is equivalent to the existence of~$\varepsilon>0$ such
that
\[
(0,\varepsilon)\subset\mathcal{A}^{\N}.
\]
\end{cor}

\begin{proof}
Continuity of the i.o.\ spectrum on $[0,a_{+})$ was shown in the
proof of~\Cref{thm:regularity-reflected-spectrum}; together with
the dimension identity in~\Cref{thm:exact-reflected-spectrum} this
proves the first assertion.

The uniform spectrum has value $\delta^{\N}(0)$ at zero and right-hand
limit $\delta^{\Z}(0)$. Put $s_{0}\coloneqq\delta^{\Z}(0)$ and $d(s)\coloneqq\partial_{q}p(s,0)$.
In the proof of~\Cref{thm:regularity-reflected-spectrum} it was
shown that $s_{0}<\delta_{0}$ and $d(s_{0})<0$ if $\alpha_{\max}<0$,
whereas $d(s_{0})\ge0$ if $\alpha_{\max}\ge0$. In the first case
also $d(\delta_{0})=\alpha_{\max}<0$, so the minimizer of $p(\delta_{0},\cdot)$
is positive, $p^{\N}_{0}(\delta_{0})=p(\delta_{0},0)=0$ by \eqref{eq:zero-speed-pressure-formulas},
and $\delta^{\N}(0)=\delta_{0}>s_{0}$. In the second case the minimizer
of $p(s_{0},\cdot)$ is non-positive, so $p^{\N}_{0}(s_{0})=p^{\Z}_{0}(s_{0})=0$
and $\delta^{\N}(0)=s_{0}$. Hence $\delta^{\N}(0)>\delta^{\Z}(0)$
exactly when $\alpha_{\max}<0$. The punctured-neighborhood characterization,
including the zero-drift case, was established in the proof of \cref{thm:regularity-reflected-spectrum}.
\end{proof}

\begin{cor}[Uniformly non-negative weighted drift]
\label{cor:nonnegative-weighted-drift} Assume 
\[
\partial_{q}p(s,0)=\int\psi\,\dd\mu_{s\phi}\ge0,\qquad0\le s\le\delta_{0}.
\]
Then, for every $\alpha\in(0,a_{+}]$, 
\[
\delta^{\N}(\alpha)=\delta^{\N}_{u}(\alpha).
\]
Consequently, $\dim_{H}E^{\N}(\alpha)=\dim_{H}E^{\N}_{u}(\alpha)=\delta^{\N}_{u}(\alpha).$
\end{cor}

\begin{proof}
By~\Cref{thm:truncated-base-pressure}, 
\[
p^{\N}_{\alpha}(s)=p^{\Z}_{\alpha}(s)\qquad(0\le s\le\delta_{0}).
\]
Both Poincaré spectra take values in $[0,\delta_{0}]$, so the dimension-relevant
zeros of the two pressure functions lie in this interval. \Cref{prop:delta-pressure-root}
therefore gives 
\[
\delta^{\N}(\alpha)=\delta^{\Z}(\alpha)=\delta^{\N}_{u}(\alpha),
\]
where the last equality is~\Cref{prop:uniform-exponent-comparison}.
The dimension identities follow from \Cref{thm:uniform-N} and~\Cref{thm:exact-reflected-spectrum}.
\end{proof}

\section{Examples}

\label{sec:examples} 

We record several explicit examples illustrating the two intrinsic
spectra of the $\N$-extension, $\delta^{\N}_{u}(\alpha)=\dim_{H}E^{\N}_{u}(\alpha)$
and $\delta^{\N}(\alpha)=\dim_{H}E^{\N}(\alpha)$. For $\alpha>0$
the uniform exponent is computed through~\Cref{prop:uniform-exponent-comparison}
as $\delta^{\N}_{u}(\alpha)=\delta^{\Z}(\alpha)$. The i.o.\ spectrum
is given by \Cref{thm:exact-reflected-spectrum}. In every example
the step set contains both signs and generates $\Z$ as an additive
group, so the standing assumptions above are satisfied; in particular,
both skew products are topologically transitive. Numerical values
are stated to the displayed precision. All inequalities, zero counts
and sign patterns asserted below that are not verified analytically
have been checked with controlled numerical precision (interval arithmetic
with outward rounding); we do not go into the details of these computations.

In the examples, $a_{+}=\max_{i\in I}\psi(i)$ is the maximal escape
rate, not the maximal-dimension drift $\alpha_{\max}$. If $\mathcal{A}^{\N}\ne\varnothing$,
\Cref{rem:no-gain-near-max-speed} gives
\[
a_{+}-\sup\mathcal{A}^{\N}>0,\qquad\delta^{\N}=\delta^{\N}_{u}\quad\text{on }[\sup\mathcal{A}^{\N},a_{+}].
\]
The width $a_{+}-\sup\mathcal{A}^{\N}$ of this terminal agreement
interval is recorded in each example; it is a horizontal distance
and must not be confused with the vertical difference $\delta^{\N}(\alpha)-\delta^{\N}_{u}(\alpha)$
plotted below. When $\mathcal{A}^{\N}=\varnothing$, the spectra agree
on the entire attainable interval $[0,a_{+}]$.

\subsection{Constant-slope examples}
\begin{figure}[t]
\centering \begin{tikzpicture}
\begin{groupplot}[
    group style={group size=2 by 1, horizontal sep=2cm},
    width=0.44\textwidth,
    height=0.34\textheight,
    xmin=0, xmax=1,
    ymin=0, ymax=1.05,
    axis lines=left,
    samples=220,
    domain=0:1,
    xtick={0,0.333333333,1},
    xticklabels={$0$,$\frac13$,$a_+=1$},
    ytick={0,0.5,1},
    xlabel={$\alpha$},
    ylabel={dimension},
    tick label style={font=\small},
    label style={font=\small},
    title style={font=\small},
]
\pgfmathdeclarefunction{HH}{1}{%
  \pgfmathparse{(#1<=0 || #1>=1) ? 0 : -(#1)*ln(#1) - (1-#1)*ln(1-#1)}%
}
\nextgroupplot[
    title={$\alpha_{\max}<0$ : $(-1,-1,+1)$},
    legend style={draw=none, fill=none, font=\small, at={(0.03,0.03)}, anchor=south west},
]
\addplot[black,dashed,thick]
{(HH((1+x)/2) + ((1-x)/2)*ln(2))/ln(3)};
\addlegendentry{$\dim_H E_{u}^{\N}(\alpha)$}
\addplot[black,solid,thick,domain=0:1/3]
{1 - x*ln(2)/ln(3)};
\addplot[black,solid,thick,domain=1/3:1]
{(HH((1+x)/2) + ((1-x)/2)*ln(2))/ln(3)};
\addlegendentry{$\dim_H E^{\N}(\alpha)$}
\addplot[black,dotted] coordinates {(1/3,0) (1/3,1)};

\nextgroupplot[
    title={$\alpha_{\max}>0$ : $(-1,+1,+1)$},
    legend style={draw=none, fill=none, font=\small, at={(0.03,0.03)}, anchor=south west},
]
\addplot[black,dashed,thick]
{(HH((1+x)/2) + ((1+x)/2)*ln(2))/ln(3)};
\addlegendentry{$\dim_H E_{u}^{\N}(\alpha)$}
\addplot[black,solid,thick]
{(HH((1+x)/2) + ((1+x)/2)*ln(2))/ln(3)};
\addlegendentry{$\dim_H E^{\N}(\alpha)$}
\addplot[black,dotted] coordinates {(1/3,0) (1/3,1)};
\node[font=\small] at (axis description cs:0.53,0.40) {$\mathcal A^{\N}=\varnothing$};
\end{groupplot}\end{tikzpicture} \caption{Two ternary constant-slope examples showing the i.o.\ and uniform
Hausdorff dimension spectra of the $\N$-extension. For $\alpha>0$
the latter is computed using $\dim_{H}E^{\N}_{u}(\alpha)=\dim_{H}E^{\Z}(\alpha)$.
In the left panel $\dim_{H}E^{\N}_{u}(0)=\dim_{H}E^{\N}(0)$, so the
uniform spectrum has a jump at zero. Left: negative maximal drift,
where the i.o.\ spectrum of the $\N$-extension develops an additional
linear branch. Right: positive maximal drift, where the i.o.\ and
uniform spectra of the $\N$-extension coincide. Both panels have
$a_{+}=1$. On the left, $\mathcal{A}^{\N}=(0,1/3)$ and the terminal
agreement interval $[1/3,1]$. On the right, $\mathcal{A}^{\N}=\varnothing$
and the spectra agree on all of $[0,1]$.}
\label{fig:constant-slope-figures}
\end{figure}
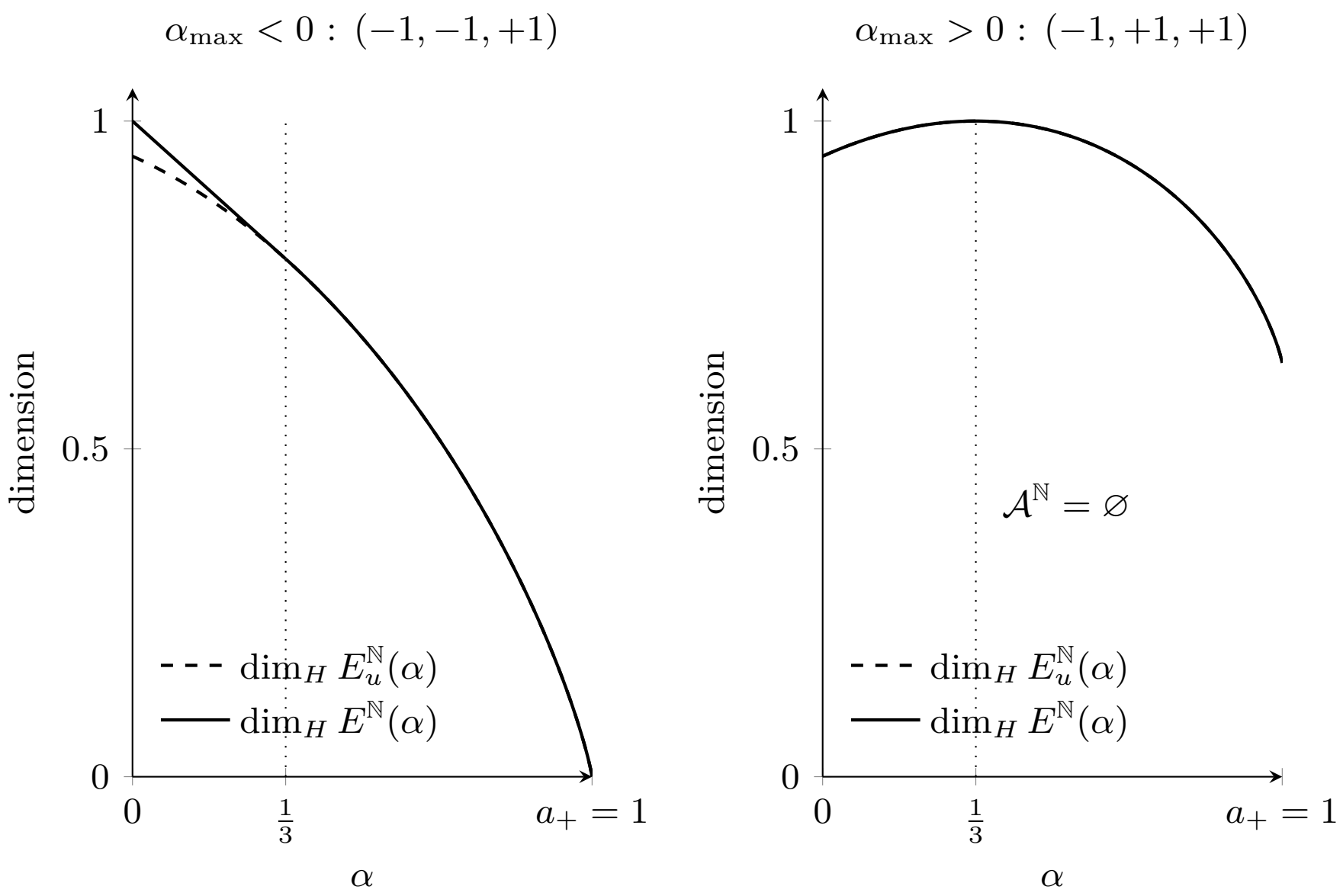

\begin{example}[A simple constant-slope family]
\label{ex:constant-slope-family} Assume that 
\[
\phi\equiv-\chi,\qquad\chi>0,
\]
and write 
\[
M\coloneqq|I|,\qquad a_{i}\coloneqq\psi(i),\qquad\Lambda(q)\coloneqq\log\sum_{i\in I}\e^{qa_{i}}.
\]
The maximal escape rate is $a_{+}=\max_{i\in I}a_{i}$. Then 
\[
p(s,q)\coloneqq\cP^{0}(s\phi+q\psi)=-s\chi+\Lambda(q),
\]
and therefore 
\[
\partial_{q}p(s,0)=\Lambda'(0)=\frac{1}{M}\sum_{i\in I}a_{i}=\alpha_{\max}
\]
is independent of $s$. Moreover, 
\[
\delta_{0}=\frac{\log M}{\chi}.
\]

For realization by disjoint full branches in $(0,1)$, the common
length $\e^{-\chi}$ satisfies $M\e^{-\chi}\le1$, equivalently $\chi\ge\log M$.

There are precisely two cases.
\begin{enumerate}
\item If $\alpha_{\max}\ge0$, then 
\[
\delta^{\N}(\alpha)=\delta^{\N}_{u}(\alpha)\qquad(0\le\alpha\le a_{+}).
\]
Consequently, 
\[
\dim_{H}E^{\N}(\alpha)=\dim_{H}E^{\N}_{u}(\alpha)=\delta^{\N}_{u}(\alpha)
\]
for every $\alpha\in(0,a_{+}]$.
\item If $\alpha_{\max}<0$, let $q^{*}>0$ be the unique nonzero solution
of 
\[
\Lambda(q^{*})=\Lambda(0)=\log M.
\]
Then 
\[
\delta^{\N}(\alpha)=\begin{cases}
{\displaystyle \delta_{0}-\frac{q^{*}}{\chi}\alpha,} & 0\le\alpha\le\Lambda'(q^{*}),\\[2mm]
\delta^{\N}_{u}(\alpha), & \Lambda'(q^{*})\le\alpha\le a_{+}.
\end{cases}
\]
In particular, the i.o.\ spectrum of the $\N$-extension has a linear
branch from $\alpha=0$ to $\alpha=\Lambda'(q^{*})$ and then joins
the uniform spectrum $\delta^{\N}_{u}$ tangentially.
\end{enumerate}
Thus $\mathcal{A}^{\N}=\varnothing$ in case (1). In case (2),
\[
\mathcal{A}^{\N}=(0,\Lambda'(q^{*})),\qquad a_{+}-\sup\mathcal{A}^{\N}=\frac{\sum_{i\in I}(a_{+}-a_{i})e^{q^{*}a_{i}}}{\sum_{i\in I}\e^{q^{*}a_{i}}}>0.
\]
Indeed, strict convexity of $\Lambda$ makes the contact branch strictly
larger before the contact escape rate $\Lambda'(q^{*})$, and the
probability vector $\bigl(\e^{q^{*}a_{i}}/\sum_{j\in I}\e^{q^{*}a_{j}}\bigr)_{i\in I}$
gives positive weight to a step smaller than $a_{+}$. Hence $[\Lambda'(q^{*}),a_{+}]$
is a non-degenerate terminal interval of equality.

For ternary examples with three branches of length $1/3$ see~\Cref{fig:constant-slope-figures}:
the left panel shows case (ii) with $\psi=(-1,-1,+1)$, where $\alpha_{\max}=-1/3$,
$q^{*}=\log2$ and $\mathcal{A}^{\N}=(0,1/3)$; the right panel shows
case (i) with $\psi=(-1,+1,+1)$, where $\alpha_{\max}=1/3$ and $\mathcal{A}^{\N}=\varnothing$.
\end{example}

\begin{proof}
If $\alpha_{\max}\ge0$, then \Cref{cor:nonnegative-weighted-drift}
gives $\delta^{\N}(\alpha)=\delta^{\N}_{u}(\alpha)$ for $\alpha\in(0,a_{+}]$,
and at $\alpha=0$ this follows from~\eqref{eq:zero-speed-pressure-formulas}.

Suppose now that $\alpha_{\max}<0$. Since 
\[
p(s,q^{*})=p(s,0)
\]
is equivalent to 
\[
\Lambda(q^{*})=\Lambda(0),
\]
the contact point $q^{*}$ is independent of $s$. The general contact
formula therefore gives 
\[
p^{\N}_{\alpha}(s)=\begin{cases}
\log M-s\chi-q^{*}\alpha, & 0\le\alpha\le\Lambda'(q^{*}),\\[1mm]
p^{\Z}_{\alpha}(s), & \Lambda'(q^{*})\le\alpha\le a_{+}.
\end{cases}
\]
Since $\delta_{0}=\log M/\chi$, taking the zero in $s$ yields the
asserted formula for $\delta^{\N}$, and tangency follows from $\left(\delta^{\N}_{u}\right)'(\Lambda'(q^{*}))=-q^{*}/\chi$.
\end{proof}

\subsection{Contact geometry and non-concavity}
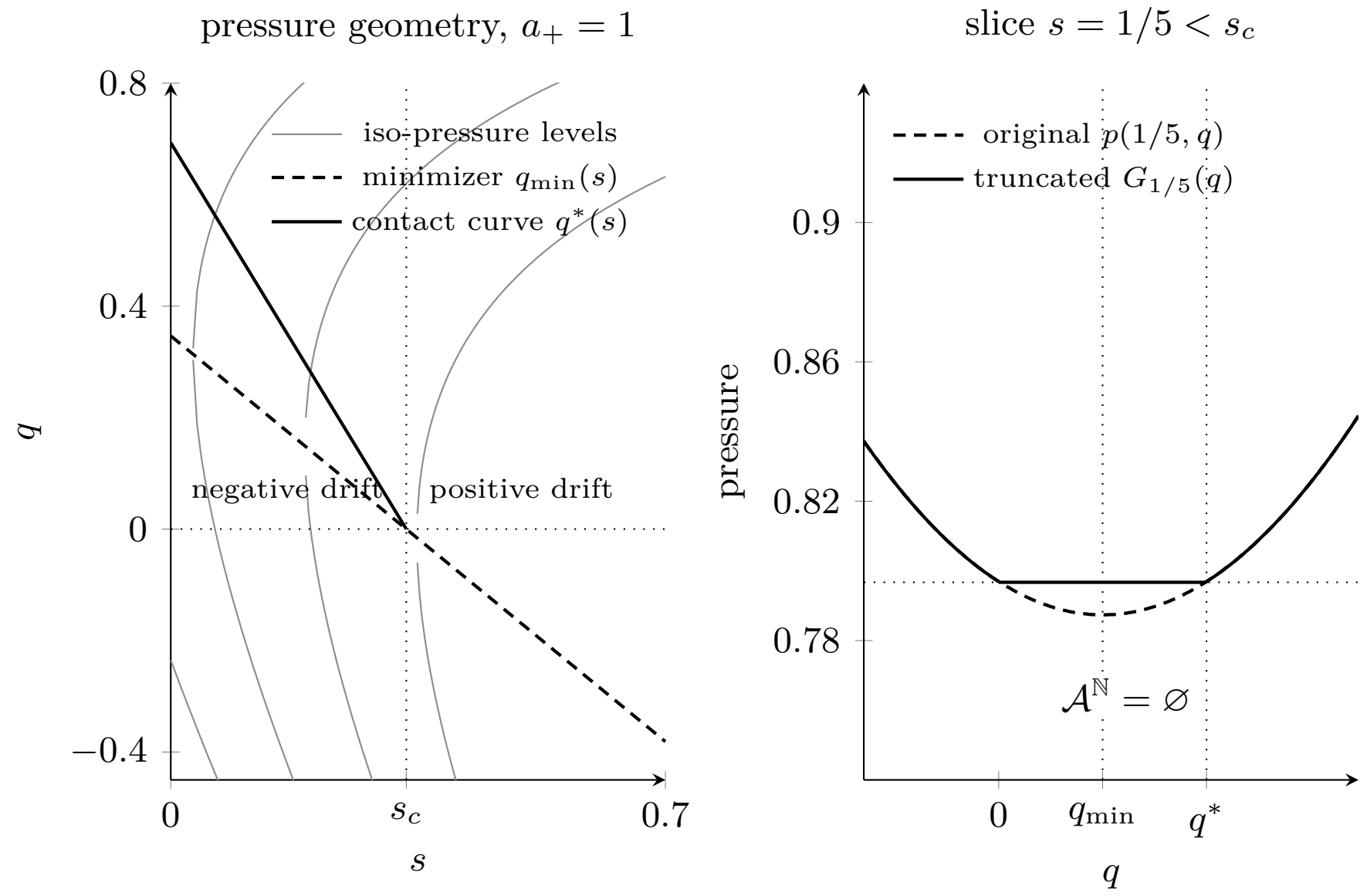
\begin{figure}[t]
\centering \begin{tikzpicture}
\begin{groupplot}[
    group style={group size=2 by 1, horizontal sep=1.7cm},
    width=0.46\textwidth,
    height=0.34\textheight,
    axis lines=left,
    tick label style={font=\small},
    label style={font=\small},
    title style={font=\small},
]

\nextgroupplot[
    xmin=0, xmax=0.7,
    ymin=-0.45, ymax=0.8,
    xlabel={$s$},
    ylabel={$q$},
    xtick={0,0.333333,0.7},
    xticklabels={$0$,$s_c$,$0.7$},
    ytick={-0.4,0,0.4,0.8},
    title={pressure geometry, $a_+=1$},
    legend pos=north east,
    legend style={draw=none,fill=none,font=\scriptsize},
]
\addplot[black!45,thin,forget plot,domain=0.349:0.7,samples=120]
{ln((exp(0.6)+sqrt(exp(1.2)-8*pow(0.08,x)))/(2*pow(0.8,x)))};
\addplot[black!45,thin,forget plot,domain=0.349:0.7,samples=120]
{ln((4*pow(0.1,x))/(exp(0.6)+sqrt(exp(1.2)-8*pow(0.08,x))))};
\addplot[black!45,thin,forget plot,domain=0.191:0.7,samples=120]
{ln((exp(0.8)+sqrt(exp(1.6)-8*pow(0.08,x)))/(2*pow(0.8,x)))};
\addplot[black!45,thin,forget plot,domain=0.191:0.7,samples=120]
{ln((4*pow(0.1,x))/(exp(0.8)+sqrt(exp(1.6)-8*pow(0.08,x))))};
\addplot[black!45,thin,forget plot,domain=0.0316:0.7,samples=120]
{ln((exp(1.0)+sqrt(exp(2.0)-8*pow(0.08,x)))/(2*pow(0.8,x)))};
\addplot[black!45,thin,forget plot,domain=0.0316:0.7,samples=120]
{ln((4*pow(0.1,x))/(exp(1.0)+sqrt(exp(2.0)-8*pow(0.08,x))))};
\addplot[black!45,thin,forget plot,domain=0:0.7,samples=120]
{ln((exp(1.2)+sqrt(exp(2.4)-8*pow(0.08,x)))/(2*pow(0.8,x)))};
\addplot[black!45,thin,forget plot,domain=0:0.7,samples=120]
{ln((4*pow(0.1,x))/(exp(1.2)+sqrt(exp(2.4)-8*pow(0.08,x))))};
\addlegendimage{black!45,thin}
\addlegendentry{iso-pressure levels}

\addplot[black,densely dashed,thick,forget plot,domain=0:0.7,samples=2]
{0.5*(ln(2)-x*ln(8))};
\addplot[black,thick,forget plot,domain=0:0.333333,samples=2]
{ln(2)-x*ln(8)};
\addlegendimage{black,densely dashed,thick}
\addlegendentry{minimizer $q_{\min}(s)$}
\addlegendimage{black,thick}
\addlegendentry{contact curve $q^*(s)$}
\addplot[black,dotted,forget plot] coordinates {(0,0) (0.7,0)};
\addplot[black,dotted,forget plot] coordinates {(0.333333,-0.45) (0.333333,0.8)};
\node[anchor=south west,font=\scriptsize] at (axis cs:0.345,0.02) {positive drift};
\node[anchor=south east,font=\scriptsize] at (axis cs:0.322,0.02) {negative drift};

\nextgroupplot[
    xmin=-0.18, xmax=0.48,
    ymin=0.74, ymax=0.94,
    xlabel={$q$},
    ylabel={pressure},
    xtick={0,0.138629,0.277259},
    xticklabels={$0$,$q_{\min}$,$q^*$},
    ytick={0.78,0.82,0.86,0.90},
    title={slice $s=1/5<s_c$},
    legend pos=north west,
    legend style={draw=none,fill=none,font=\scriptsize},
]
\addplot[black,densely dashed,thick,forget plot,domain=-0.18:0.48,samples=180]
{ln(2*pow(0.1,0.2)*exp(-x)+pow(0.8,0.2)*exp(x))};
\addplot[black,thick,forget plot,domain=-0.18:0,samples=90]
{ln(2*pow(0.1,0.2)*exp(-x)+pow(0.8,0.2)*exp(x))};
\addplot[black,thick,forget plot] coordinates {(0,0.7967263455) (0.2772588722,0.7967263455)};
\addplot[black,thick,forget plot,domain=0.2772588722:0.48,samples=90]
{ln(2*pow(0.1,0.2)*exp(-x)+pow(0.8,0.2)*exp(x))};
\addlegendimage{black,densely dashed,thick}
\addlegendentry{original $p(1/5,q)$}
\addlegendimage{black,thick}
\addlegendentry{truncated $G_{1/5}(q)$}
\addplot[black,dotted,forget plot] coordinates {(0.1386294361,0.74) (0.1386294361,0.94)};
\addplot[black,dotted,forget plot] coordinates {(0.2772588722,0.74) (0.2772588722,0.94)};
\addplot[black,dotted,forget plot] coordinates {(-0.18,0.7967263455) (0.48,0.7967263455)};
\node[font=\small,fill=white,inner sep=1pt] at (axis description cs:0.53,0.12) {$\mathcal A^{\N}=\varnothing$};
\end{groupplot}\end{tikzpicture} \caption{Geometry behind \Cref{thm:truncated-base-pressure} for the example
in \Cref{ex:loc-contact-geometry-1}. Left: iso-pressure curves of
$(s,q)\protect\mapsto p(s,q)$ together with the minimizer $q_{\min}(s)$
and the second contact point $q^{*}(s)$. The positive part of the
contact curve exists precisely for $\partial_{q}p(s,0)<0$ and reaches
zero at $s_{c}=1/3$. Right: for the slice $s=1/5$, truncation replaces
the dip of the dashed original pressure $p(1/5,q)$ between $0$ and
$q^{*}$ by the horizontal part of the solid truncated pressure $G_{1/5}(q)$
at level $p(1/5,0)$; its Legendre transform therefore has the linear
contact branch in \eqref{eq:contact-formula}. Here $a_{+}=1$ and
$\mathcal{A}^{\N}=\varnothing$: the pressure contact range displayed
here does not produce strict gain at the dimension zeros, as proved
in the example.}
\label{fig:loc-contact-geometry-1}
\end{figure}

\begin{example}[Geometry of the contact formula]
\label{ex:loc-contact-geometry-1} The two cases in \Cref{thm:truncated-base-pressure}
can be seen directly in the $(s,q)$-pressure plane. Consider three
linear full branches with 
\[
r_{1}=r_{2}=\frac{1}{10},\qquad r_{3}=\frac{4}{5},\qquad\psi(1)=\psi(2)=-1,\qquad\psi(3)=+1.
\]
The maximal escape rate is $a_{+}=1$. Then 
\[
p(s,q)=\log\left(2\,10^{-s}\e^{-q}+\left(\frac{4}{5}\right)^{s}\e^{q}\right).
\]
The drift $d(s)=\partial_{q}p(s,0)$ changes sign at 
\[
s_{c}=\frac{\log2}{\log8}=\frac{1}{3},\qquad\partial_{q}p(s,0)<0\ \text{for }s<s_{c},\qquad\partial_{q}p(s,0)>0\ \text{for }s>s_{c}.
\]
The unconstrained minimizer is 
\[
q_{\min}(s)=\frac{1}{2}\bigl(\log2-s\log8\bigr),
\]
and, for $0\le s<s_{c}$, the second intersection of the level through
$q=0$ is 
\[
q^{*}(s)=\log2-s\log8=2q_{\min}(s).
\]
Thus $p(s,q^{*}(s))=p(s,0)$, and the truncation $G_{s}(q)=\max\{p(s,0),p(s,q)\}$
is flat on $[0,q^{*}(s)]$. At $s=s_{c}$ both $q_{\min}(s)$ and
$q^{*}(s)$ meet the line $q=0$; for $s>s_{c}$ the minimizer lies
in $q<0$, so the truncation no longer changes the Legendre transform
for $\alpha>0$. This is the transition between the two alternatives
in \Cref{thm:truncated-base-pressure}. Both the iso-pressure geometry
and the truncation for a fixed $s$ are shown in~\cref{fig:loc-contact-geometry-1}.

In this example the pressure truncation does not create a dimension
gap: $\mathcal{A}^{\N}=\varnothing$, so the spectra agree on all
of $[0,a_{+}]=[0,1]$. To see this directly, for $0\le s<s_{c}$ the
formulas above give $0<q^{*}(s)\le\log2$ and $0<a^{*}(s)\le1/3$.
Since $p(s,0)$ is decreasing,
\[
H(s)=p(s,0)-q^{*}(s)a^{*}(s)\ge\log\!\left(2(4/5)^{1/3}\right)-\frac{\log2}{3}=\frac{1}{3}\log\frac{16}{5}>0.
\]
A strict-gain zero would have $0<\alpha<a^{*}(s)$ and $p^{\N}_{\alpha}(s)=p(s,0)-q^{*}(s)\alpha>H(s)>0$,
a contradiction. For $s\ge s_{c}$, $d(s)\ge0$ gives $a^{*}(s)\le0$,
so strict gain is excluded by~\Cref{cor:strict-reflection-gain}.
The contact segments in the pressure diagram therefore do not meet
a dimension zero in their strict-contact range.
\end{example}

\begin{example}[A non-concave i.o.\ dimension spectrum for the $\N$-extension]
\label{ex:nonconcave-reflected} Consider two linear full branches
with geometric lengths $r_{-}\coloneqq9/10$, $r_{+}\coloneqq1/10$,
and steps $\psi(-)\coloneqq-1$, $\psi(+)\coloneqq+1$. Thus $a_{+}=1$.
Then $\delta_{0}=1$ and $\alpha_{\max}=-4/5.$ At zero, $\delta^{\N}_{u}(0)=\delta^{\N}(0)=1.$
The uniform branch for $\alpha>0$ has a smaller right-hand limit
at $0$. For $0<\alpha\leq1$, writing $\theta=(1+\alpha)/2$, the
uniform exponent is 
\[
\delta^{\N}_{u}(\alpha)=\delta^{\Z}(\alpha)=\frac{\mathscr{H}(\theta)}{-(1-\theta)\log(9/10)-\theta\log(1/10)},
\]
where $\mathscr{H}(\theta)\coloneqq-\theta\log\theta-(1-\theta)\log(1-\theta).$
The locally constant pressure is 
\[
p(s,q)=\log\left(\left(\frac{9}{10}\right)^{s}\e^{-q}+\left(\frac{1}{10}\right)^{s}\e^{q}\right).
\]
For every $s>0$ the drift $d(s)=\partial_{q}p(s,0)$ is negative,
and therefore 
\[
p^{\N}_{0}(s)=p(s,0)=:f(s)=\log\left(\left(\frac{9}{10}\right)^{s}+\left(\frac{1}{10}\right)^{s}\right).
\]
The equation $p(s,q)=p(s,0)$ has, besides $q=0$, the second root
$q^{*}(s)=s\log9,$ and hence the contact escape rate is $\partial_{q}p(s,q^{*}(s))=\left((9/10)^{s}-(1/10)^{s}\right)/\left((9/10)^{s}+(1/10)^{s}\right).$
Thus \Cref{thm:truncated-base-pressure} gives 
\[
p^{\N}_{\alpha}(s)=\begin{cases}
f(s)-\alpha s\log9, & 0\leq\alpha\leq\partial_{q}p(s,q^{*}(s)),\\[1mm]
p^{\Z}_{\alpha}(s), & \partial_{q}p(s,q^{*}(s))\leq\alpha\leq1.
\end{cases}
\]
Let $s_{*}$ be the unique solution of $f(s_{*})=s_{*}\log9\,\partial_{q}p(s_{*},q^{*}(s_{*})),$
and put $\alpha^{*}\coloneqq\partial_{q}p(s_{*},q^{*}(s_{*}))$. Numerically,
$s_{*}\approx0.3761469$, $\alpha^{*}\approx0.3912199$. Consequently
\[
\delta^{\N}(\alpha)=\begin{cases}
s(\alpha), & 0\leq\alpha\leq\alpha^{*},\\[1mm]
\delta^{\N}_{u}(\alpha), & \alpha^{*}\leq\alpha\leq1,
\end{cases}
\]
where $s(\alpha)\in[s_{*},1]$ is the unique solution of $f(s)=\alpha s\log9.$
Moreover, $H'(s^{*})\neq0$, where $H(s):=f(s)-s\log9\,a^{*}(s)$.
Hence \Cref{cor:contact-multiplicity} shows that $\alpha^{*}$ is
a second-order phase transition: $\delta^{\N}$ is $C^{1}$ but not
$C^{2}$ at $\alpha^{*}$. Note that $\alpha\mapsto\delta^{\N}(\alpha)$
is not concave. Indeed, $s(0)=1$ and, with 
\[
d\coloneqq f'(1)=\frac{9}{10}\log\frac{9}{10}+\frac{1}{10}\log\frac{1}{10}<0,
\]
implicit differentiation of $\alpha\mapsto s\left(\alpha\right)$
gives $(\delta^{\N})'(0+)=s'(0)=\log\left(9\right)/d<0$. Moreover
$f''(1)=9/100\cdot(\log9)^{2}>0$ and 
\[
(\delta^{\N})''(0+)=-\frac{f''(1)s'(0)^{2}-2\log9\,s'(0)}{d}>0.
\]
Hence $\delta^{\N}$ is strictly convex in a right neighbourhood of
$0$ and therefore cannot be concave on $[0,1]$. Both spectra are
plotted in~\cref{fig:nonconcave-reflected}. Here $\mathcal{A}^{\N}=(0,\alpha^{*})$
and $a_{+}-\sup\mathcal{A}^{\N}=1-\alpha^{*}\approx0.6088$; see~\Cref{fig:nonconcave-reflected}.
\end{example}

\subsection{Strict gain despite positive maximal-dimension drift}
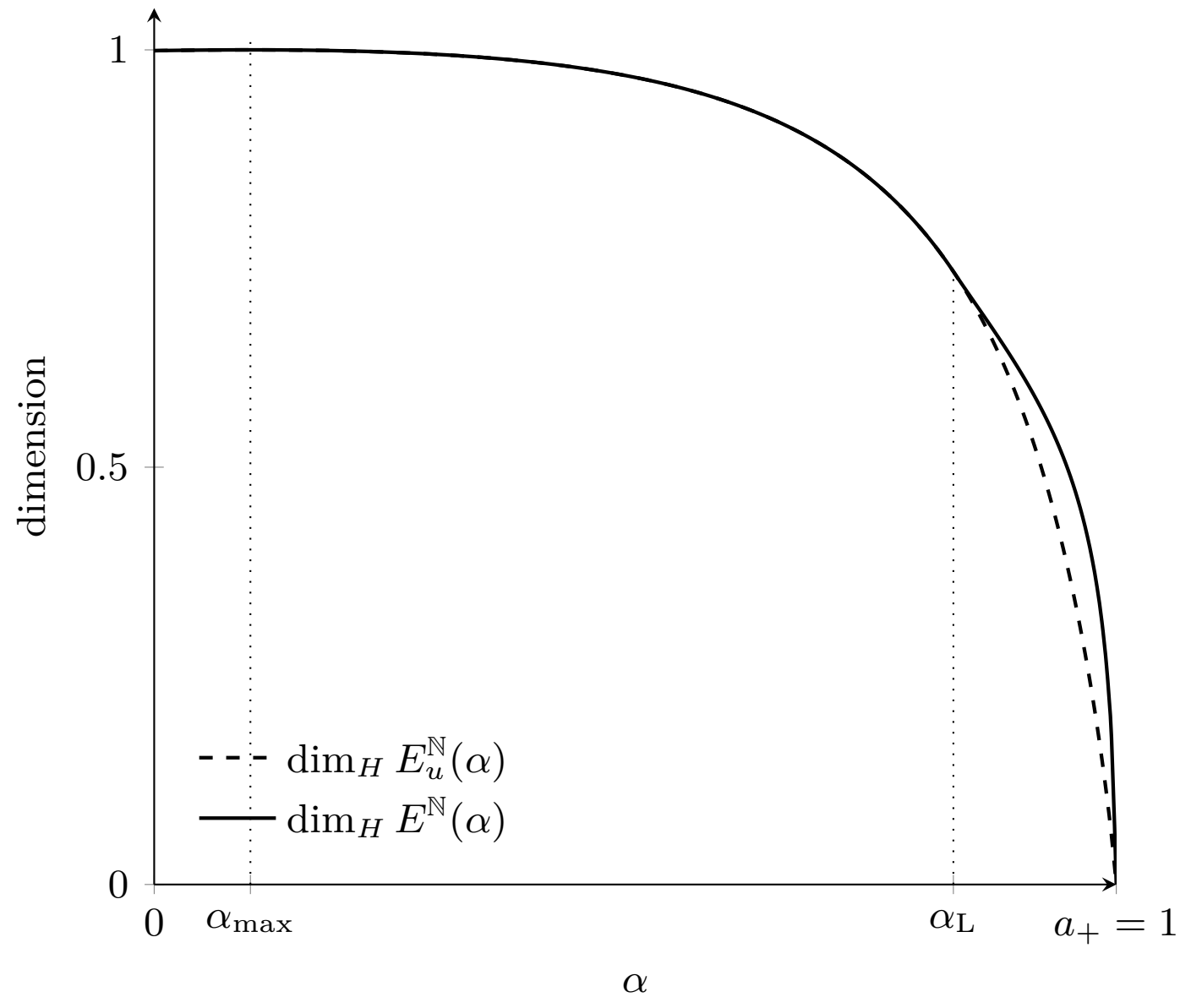
\begin{figure}[t]
\centering \begin{tikzpicture}
\begin{axis}[
    width=0.72\textwidth,
    height=0.38\textheight,
    xmin=0, xmax=1,
    ymin=0, ymax=1.05,
    axis lines=left,
    xlabel={$\alpha$},
    ylabel={dimension},
    xtick={0,0.1,0.8306633121292714,1},
    xticklabels={$0$,$\alpha_\mathrm{max} $,$\alpha_{\mathrm{L}}$,$a_+=1$},
    ytick={0,0.5,1},
    tick label style={font=\small},
    label style={font=\small},
    legend style={draw=none, fill=none, font=\small, at={(0.03,0.03)}, anchor=south west},
]
\addplot[black,dashed,thick,smooth] coordinates {
(0.000,0.999028) (0.005,0.999118) (0.010,0.999204) (0.015,0.999286) (0.020,0.999363) (0.025,0.999437) (0.030,0.999507) (0.035,0.999572)
(0.040,0.999633) (0.045,0.999690) (0.050,0.999742) (0.055,0.999790) (0.060,0.999833) (0.065,0.999871) (0.070,0.999905) (0.075,0.999933)
(0.080,0.999957) (0.085,0.999976) (0.090,0.999989) (0.095,0.999997) (0.100,1.000000) (0.105,0.999997) (0.110,0.999989) (0.115,0.999975)
(0.120,0.999955) (0.125,0.999929) (0.130,0.999897) (0.135,0.999859) (0.140,0.999815) (0.145,0.999764) (0.150,0.999707) (0.155,0.999643)
(0.160,0.999572) (0.165,0.999495) (0.170,0.999410) (0.175,0.999318) (0.180,0.999218) (0.185,0.999112) (0.190,0.998997) (0.195,0.998875)
(0.200,0.998744) (0.205,0.998606) (0.210,0.998459) (0.215,0.998304) (0.220,0.998140) (0.225,0.997967) (0.230,0.997786) (0.235,0.997595)
(0.240,0.997395) (0.245,0.997185) (0.250,0.996965) (0.255,0.996735) (0.260,0.996495) (0.265,0.996245) (0.270,0.995984) (0.275,0.995712)
(0.280,0.995429) (0.285,0.995135) (0.290,0.994829) (0.295,0.994511) (0.300,0.994181) (0.305,0.993839) (0.310,0.993484) (0.315,0.993116)
(0.320,0.992735) (0.325,0.992341) (0.330,0.991932) (0.335,0.991510) (0.340,0.991073) (0.345,0.990621) (0.350,0.990155) (0.355,0.989672)
(0.360,0.989174) (0.365,0.988660) (0.370,0.988130) (0.375,0.987582) (0.380,0.987017) (0.385,0.986434) (0.390,0.985834) (0.395,0.985214)
(0.400,0.984576) (0.405,0.983918) (0.410,0.983240) (0.415,0.982542) (0.420,0.981823) (0.425,0.981083) (0.430,0.980320) (0.435,0.979535)
(0.440,0.978728) (0.445,0.977896) (0.450,0.977040) (0.455,0.976160) (0.460,0.975254) (0.465,0.974321) (0.470,0.973363) (0.475,0.972376)
(0.480,0.971362) (0.485,0.970318) (0.490,0.969245) (0.495,0.968142) (0.500,0.967007) (0.505,0.965840) (0.510,0.964640) (0.515,0.963407)
(0.520,0.962138) (0.525,0.960834) (0.530,0.959493) (0.535,0.958115) (0.540,0.956698) (0.545,0.955241) (0.550,0.953743) (0.555,0.952202)
(0.560,0.950619) (0.565,0.948991) (0.570,0.947317) (0.575,0.945595) (0.580,0.943825) (0.585,0.942005) (0.590,0.940134) (0.595,0.938209)
(0.600,0.936229) (0.605,0.934193) (0.610,0.932099) (0.615,0.929945) (0.620,0.927729) (0.625,0.925449) (0.630,0.923104) (0.635,0.920690)
(0.640,0.918206) (0.645,0.915650) (0.650,0.913019) (0.655,0.910311) (0.660,0.907523) (0.665,0.904652) (0.670,0.901696) (0.675,0.898651)
(0.680,0.895515) (0.685,0.892284) (0.690,0.888955) (0.695,0.885525) (0.700,0.881990) (0.705,0.878346) (0.710,0.874589) (0.715,0.870714)
(0.720,0.866718) (0.725,0.862596) (0.730,0.858343) (0.735,0.853954) (0.740,0.849423) (0.745,0.844745) (0.750,0.839915) (0.755,0.834925)
(0.760,0.829770) (0.765,0.824442) (0.770,0.818934) (0.775,0.813240) (0.780,0.807349) (0.785,0.801255) (0.790,0.794948) (0.795,0.788419)
(0.800,0.781657) (0.805,0.774652) (0.810,0.767393) (0.815,0.759867) (0.820,0.752063) (0.825,0.743965) (0.830,0.735560) (0.835,0.726833)
(0.840,0.717767) (0.845,0.708344) (0.850,0.698546) (0.855,0.688353) (0.860,0.677744) (0.865,0.666694) (0.870,0.655180) (0.875,0.643175)
(0.880,0.630650) (0.885,0.617573) (0.890,0.603912) (0.895,0.589629) (0.900,0.574686) (0.905,0.559038) (0.910,0.542639) (0.915,0.525436)
(0.920,0.507373) (0.925,0.488387) (0.930,0.468408) (0.935,0.447359) (0.940,0.425154) (0.945,0.401694) (0.950,0.376869) (0.955,0.350552)
(0.960,0.322596) (0.965,0.292827) (0.970,0.261039) (0.975,0.226973) (0.980,0.190301) (0.985,0.150572) (0.990,0.107103) (0.995,0.058626)
(1.000,0.000000)
};
\addlegendentry{$\dim_H E_{u}^{\N}(\alpha)$}
\addplot[black,solid,thick,smooth] coordinates {
(0.000,0.999028) (0.005,0.999118) (0.010,0.999204) (0.015,0.999286) (0.020,0.999363) (0.025,0.999437) (0.030,0.999507) (0.035,0.999572)
(0.040,0.999633) (0.045,0.999690) (0.050,0.999742) (0.055,0.999790) (0.060,0.999833) (0.065,0.999871) (0.070,0.999905) (0.075,0.999933)
(0.080,0.999957) (0.085,0.999976) (0.090,0.999989) (0.095,0.999997) (0.100,1.000000) (0.105,0.999997) (0.110,0.999989) (0.115,0.999975)
(0.120,0.999955) (0.125,0.999929) (0.130,0.999897) (0.135,0.999859) (0.140,0.999815) (0.145,0.999764) (0.150,0.999707) (0.155,0.999643)
(0.160,0.999572) (0.165,0.999495) (0.170,0.999410) (0.175,0.999318) (0.180,0.999218) (0.185,0.999112) (0.190,0.998997) (0.195,0.998875)
(0.200,0.998744) (0.205,0.998606) (0.210,0.998459) (0.215,0.998304) (0.220,0.998140) (0.225,0.997967) (0.230,0.997786) (0.235,0.997595)
(0.240,0.997395) (0.245,0.997185) (0.250,0.996965) (0.255,0.996735) (0.260,0.996495) (0.265,0.996245) (0.270,0.995984) (0.275,0.995712)
(0.280,0.995429) (0.285,0.995135) (0.290,0.994829) (0.295,0.994511) (0.300,0.994181) (0.305,0.993839) (0.310,0.993484) (0.315,0.993116)
(0.320,0.992735) (0.325,0.992341) (0.330,0.991932) (0.335,0.991510) (0.340,0.991073) (0.345,0.990621) (0.350,0.990155) (0.355,0.989672)
(0.360,0.989174) (0.365,0.988660) (0.370,0.988130) (0.375,0.987582) (0.380,0.987017) (0.385,0.986434) (0.390,0.985834) (0.395,0.985214)
(0.400,0.984576) (0.405,0.983918) (0.410,0.983240) (0.415,0.982542) (0.420,0.981823) (0.425,0.981083) (0.430,0.980320) (0.435,0.979535)
(0.440,0.978728) (0.445,0.977896) (0.450,0.977040) (0.455,0.976160) (0.460,0.975254) (0.465,0.974321) (0.470,0.973363) (0.475,0.972376)
(0.480,0.971362) (0.485,0.970318) (0.490,0.969245) (0.495,0.968142) (0.500,0.967007) (0.505,0.965840) (0.510,0.964640) (0.515,0.963407)
(0.520,0.962138) (0.525,0.960834) (0.530,0.959493) (0.535,0.958115) (0.540,0.956698) (0.545,0.955241) (0.550,0.953743) (0.555,0.952202)
(0.560,0.950619) (0.565,0.948991) (0.570,0.947317) (0.575,0.945595) (0.580,0.943825) (0.585,0.942005) (0.590,0.940134) (0.595,0.938209)
(0.600,0.936229) (0.605,0.934193) (0.610,0.932099) (0.615,0.929945) (0.620,0.927729) (0.625,0.925449) (0.630,0.923104) (0.635,0.920690)
(0.640,0.918206) (0.645,0.915650) (0.650,0.913019) (0.655,0.910311) (0.660,0.907523) (0.665,0.904652) (0.670,0.901696) (0.675,0.898651)
(0.680,0.895515) (0.685,0.892284) (0.690,0.888955) (0.695,0.885525) (0.700,0.881990) (0.705,0.878346) (0.710,0.874589) (0.715,0.870714)
(0.720,0.866718) (0.725,0.862596) (0.730,0.858343) (0.735,0.853954) (0.740,0.849423) (0.745,0.844745) (0.750,0.839915) (0.755,0.834925)
(0.760,0.829770) (0.765,0.824442) (0.770,0.818934) (0.775,0.813240) (0.780,0.807349) (0.785,0.801255) (0.790,0.794948) (0.795,0.788419)
(0.800,0.781657) (0.805,0.774652) (0.810,0.767393) (0.815,0.759867) (0.820,0.752063) (0.825,0.743965) (0.830,0.735560) (0.835,0.727008)
(0.840,0.718553) (0.845,0.710171) (0.850,0.701830) (0.855,0.693498) (0.860,0.685154) (0.865,0.676776) (0.870,0.668336) (0.875,0.659810)
(0.880,0.651179) (0.885,0.642412) (0.890,0.633481) (0.895,0.624359) (0.900,0.615005) (0.905,0.605392) (0.910,0.595474) (0.915,0.585191)
(0.920,0.574493) (0.925,0.563314) (0.930,0.551557) (0.935,0.539131) (0.940,0.525905) (0.945,0.511725) (0.950,0.496390) (0.955,0.479638)
(0.960,0.461133) (0.965,0.440396) (0.970,0.416765) (0.975,0.389248) (0.980,0.356243) (0.985,0.314969) (0.990,0.259752) (0.995,0.176111)
(1.000,0.000000)
};
\addlegendentry{$\dim_H E^{\N}(\alpha)$}
\addplot[black,dotted] coordinates {(0.1,0) (0.1,1.02)};
\addplot[black,dotted] coordinates {(0.8306633121292714,0) (0.8306633121292714,0.7344213794442977)};
\end{axis}\end{tikzpicture} \caption{The positive-drift system of \cref{ex:loc-positive-counter}, with
$\alpha_{\max}=1/10$ and maximal escape rate $a_{+}=1$. The i.o.\
and uniform spectra agree on $[0,\alpha_{\mathrm{L}}]\cup[\alpha_{\mathrm{R}},1]$
and differ exactly on $\mathcal{A}^{\N}=(\alpha_{\mathrm{L}},\alpha_{\mathrm{R}})$.
Here $\alpha_{\mathrm{L}}\approx0.83066331$ and $1-\alpha_{\mathrm{R}}\approx4.73605\times10^{-23}$.
The left contact is marked by the dotted line. The terminal agreement
interval has positive horizontal width but is not resolved at this
scale; both spectra equal zero at $a_{+}=1$.}
\label{fig:loc-positive-counter}
\end{figure}

\begin{example}[Positive drift with exactly one strict-gain interval]
\label{ex:loc-positive-counter}Consider three full linear branches
with lengths $r_{1}=r_{2}=3/340$, $r_{3}=167/170$, $\psi(1)=\psi(2)=-50$,
$\quad\psi(3)=1.$ Then $a_{+}=1$, $\delta_{0}=1$ and $\alpha_{\max}=1/10>0$.
Despite this positive maximal-dimension drift, the strict-gain set
is exactly one nonempty interval:
\begin{equation}
\mathcal{A}^{\N}=(\alpha_{\mathrm{L}},\alpha_{\mathrm{R}}),\qquad\alpha_{\mathrm{L}}\approx0.83066331.\label{eq:positive-exact-gain-interval}
\end{equation}
The right endpoint is strictly below the maximal escape rate, with
$1-\alpha_{\mathrm{R}}\approx4.73605\times10^{-23}.$ Both boundary
points are simple contacts, so the i.o.\ spectrum is $C^{1}$ but
not $C^{2}$ at each of them.

For $0<\alpha\leq1$, put $\theta=(1-\alpha)/51$. The uniform spectrum
is given, with $\mathscr{H}$ as above, by
\[
\begin{aligned}\delta^{\N}_{u}(\alpha) & =\frac{\mathscr{H}(\theta)+\theta\log2}{-\theta\log(3/340)-(1-\theta)\log(167/170)}.\end{aligned}
\]
The i.o.\ spectrum is obtained from the contact formula. The terminal
agreement interval $[\alpha_{\mathrm{R}},1]$ is not resolved at the
scale of~\Cref{fig:loc-positive-counter}.

To see \cref{eq:positive-exact-gain-interval}, write $c(s)=p(s,0)$
and $H(s)=c(s)-q^{*}(s)a^{*}(s)$. The contact formula gives
\begin{equation}
\mathcal{A}^{\N}=\left\{ \frac{c(s)}{q^{*}(s)}:s\geq0,\ q^{*}(s)>0,\ c(s)>0,\ H(s)<0\right\} .\label{eq:gain-contact-parametrisation}
\end{equation}
Indeed, these conditions place $\alpha=c(s)/q^{*}(s)$ strictly inside
the contact branch and make its pressure zero. Conversely, for $\alpha\in\mathcal{A}^{\N}$
the contact formula gives $c(s)=\alpha q^{*}(s)>0$ and $H(s)=q^{*}(s)(\alpha-a^{*}(s))<0$,
with $s=\delta^{\N}(\alpha)$; see \cref{cor:strict-reflection-gain}.
\end{example}

\begin{proof}
Here $q^{*}(s)>0$ exactly for $0\le s<s_{c}\coloneqq\log100/\log(334/3)$,
and $c(s)>0$ on this range. With controlled numerical precision one
verifies that $H$ has exactly two zeros $s_{\mathrm{R}}<s_{\mathrm{L}}$
in $[0,s_{c})$, that both are simple, and that $H<0$ exactly on
$(s_{\mathrm{R}},s_{\mathrm{L}})$. On the active contact branch the
pressure $p^{\N}_{\alpha}(s)=c(s)-\alpha q^{*}(s)$ is strictly decreasing
in $s$ and in $\alpha$, so $s\mapsto c(s)/q^{*}(s)$ is strictly
decreasing there. Applying~\eqref{eq:gain-contact-parametrisation}
therefore gives~\eqref{eq:positive-exact-gain-interval} with $\alpha_{\mathrm{L}}=c(s_{\mathrm{L}})/q^{*}(s_{\mathrm{L}})$
and $\alpha_{\mathrm{R}}=c(s_{\mathrm{R}})/q^{*}(s_{\mathrm{R}})$,
and~\Cref{cor:contact-multiplicity} gives the regularity assertion.
\end{proof}

\begin{rem}[Three branches do not force a connected gain set]
\label{rem:three-branch-disconnected-gain}Exactly one gain interval
is a feature of the preceding parameters, not a consequence of having
three branches. Consider instead
\[
(r_{1},r_{2},r_{3})=\frac{1}{25000}(50,1,24949),\qquad(\psi(1),\psi(2),\psi(3))=(-400,-40,1).
\]
Again $a_{+}=1$ and $\delta_{0}=1$, while $\alpha_{\max}=4909/25000>0$.
For $p(s,q)=\log\sum_{i}r^{s}_{i}e^{q\psi(i)}$, use the same notation
$c(s)$, $q^{*}(s)$ and $H(s)$. Controlled-precision evaluation
gives $c(s)>0$ and $q^{*}(s)>0$ at $s=1/5,1/2,7/10$, with $H(1/5)<0$,
$H(1/2)>0$, $H(7/10)<0$. By \cref{eq:gain-contact-parametrisation},
the escape rates $\alpha_{a}=c(7/10)/q^{*}(7/10)$, $\alpha_{b}=c(1/5)/q^{*}(1/5)$
are strict-gain points with $0<\alpha_{a}<\alpha_{b}<1$. They satisfy
$\delta^{\N}(\alpha_{a})=7/10$ and $\delta^{\N}(\alpha_{b})=1/5$.
If they belonged to one component of $\mathcal{A}^{\N}$, continuity
of $\delta^{\N}$ would give a strict-gain point between them with
dimension $1/2$, forcing $H(1/2)<0$. This contradiction proves that
the three-branch system has at least two strict-gain components despite
its positive maximal-dimension drift.
\end{rem}

\subsection{Disconnected strict-gain sets}
\begin{example}[Exactly three strict-gain components]
\label{ex:three-component-gain-set} Consider ten full linear branches
with lengths 
\[
(r_{1},\ldots,r_{10})=\frac{1}{10^{6}}(2,906200,66672,446,5108,3681,13206,941,866,2878)
\]
 and integer steps 
\[
(a_{1},\ldots,a_{10})=(-85,17,-227,-54,-273,-55,-189,-92,-58,-288).
\]
In particular, $a_{+}=17$ and the attainable escape rates form $[0,17]$.
The lengths sum to one, so $\delta_{0}=1$. Branch $2$ is the only
positive-step branch and is strictly longer than every other branch.
Thus $s\mapsto\left(\sum_{i\ne2}|a_{i}|r^{s}_{i}\right)/17r^{s}_{2}$
is decreasing. At $s=1$, its numerator and denominator, multiplied
by $10^{6}$, are $20217335$ and $15405400$, respectively. Hence
$d(s)=\partial_{q}p(s,0)<0$ on $[0,1]$; in particular $\alpha_{\max}<0$
and $q^{*}(s)>0$ on $[0,1]$. The strict-gain set has exactly three
components: 
\begin{equation}
\mathcal{A}^{\N}=(0,\beta_{1})\cup(\beta_{2},\beta_{3})\cup(\beta_{4},\beta_{5}),\label{eq:three-exact-gain-components}
\end{equation}
 with approximate endpoints 
\[
\begin{array}{ccccc}
\beta_{1}\approx7.3018, & \beta_{2}\approx12.2437, & \beta_{3}\approx15.5591, & \beta_{4}\approx16.7450, & \beta_{5}\approx16.9807.\end{array}
\]
 To define them exactly, let $t_{1}<\cdots<t_{5}$ be the zeros of
$H$ in $[0,1]$. With controlled numerical precision one verifies
that these are the only zeros of $H$ in $[0,1]$, that all of them
are simple, that the sign pattern of $H$ on $[0,1]$ is $+,-,+,-,+,-$,
and that $\delta^{\Z}(0)>0.99>t_{5}$. For $s\in[0,\delta^{\Z}(0))$
let $\alpha_{\Z}(s)\in(0,a_{+}]$ be the solution of $p^{\Z}_{\alpha_{\Z}(s)}(s)=0$;
since $\alpha_{\max}<0$, the map $\delta^{\Z}$ is strictly decreasing
on $[0,a_{+}]$ and $\alpha_{\Z}$ is its inverse. As $\alpha\mapsto p^{\Z}_{\alpha}(s)$
is strictly decreasing on $[0,a_{+}]$ (its maximum is attained at
$\alpha=d(s)<0$) and $H(s)=p^{\Z}_{a^{*}(s)}(s)$, we have $H(s)<0$
if and only if $\alpha_{\Z}(s)<a^{*}(s)$, which by~\Cref{cor:strict-reflection-gain}
is equivalent to $\alpha_{\Z}(s)\in\mathcal{A}^{\N}$; moreover, $a^{*}(t)=\alpha_{\Z}(t)$
at every zero $t$ of $H$. Setting $\beta_{i}\coloneqq a^{*}(t_{6-i})=\alpha_{\Z}(t_{6-i})$,
the sign pattern of $H$ therefore proves~\eqref{eq:three-exact-gain-components}.
Since all five zeros of $H$ are simple, all five transitions are
$C^{1}$ but not $C^{2}$ by~\Cref{cor:contact-multiplicity}. These
separated components arise although the drift $d(s)$ never changes
sign on $[0,1]$; see~\Cref{fig:three-component-gain-set}.

Here $a_{+}-\sup\mathcal{A}^{\N}=17-\beta_{5}\approx0.0193$, so the
third gain component ends before the maximal escape rate $a_{+}=17$,
as the lower panel of~\Cref{fig:three-component-gain-set} shows.
\end{example}

\emph{Data availability statement:} No new data were created or analysed
in this study.

We declare that there is no conflict of interest.

\end{document}